\expandafter\chardef\csname pre amssym.def at\endcsname=\the\catcode`\@
\catcode`\@=11

\def\undefine#1{\let#1\undefined}
\def\newsymbol#1#2#3#4#5{\let\next@\relax
 \ifnum#2=\@ne\let\next@\msafam@\else
 \ifnum#2=\tw@\let\next@\msbfam@\fi\fi
 \mathchardef#1="#3\next@#4#5}
\def\mathhexbox@#1#2#3{\relax
 \ifmmode\mathpalette{}{\m@th\mathchar"#1#2#3}%
 \else\leavevmode\hbox{$\m@th\mathchar"#1#2#3$}\fi}
\def\hexnumber@#1{\ifcase#1 0\or 1\or 2\or 3\or 4\or 5\or 6\or 7\or 8\or
 9\or A\or B\or C\or D\or E\or F\fi}

\font\tenmsa=msam10
\font\sevenmsa=msam7
\font\fivemsa=msam5
\newfam\msafam
\textfont\msafam=\tenmsa
\scriptfont\msafam=\sevenmsa
\scriptscriptfont\msafam=\fivemsa
\edef\msafam@{\hexnumber@\msafam}
\mathchardef\dabar@"0\msafam@39
\def\dashrightarrow{\mathrel{\dabar@\dabar@\mathchar"0\msafam@4B}}
\def\dashleftarrow{\mathrel{\mathchar"0\msafam@4C\dabar@\dabar@}}

\def\ulcorner{\delimiter"4\msafam@70\msafam@70 }
\def\urcorner{\delimiter"5\msafam@71\msafam@71 }
\def\llcorner{\delimiter"4\msafam@78\msafam@78 }
\def\lrcorner{\delimiter"5\msafam@79\msafam@79 }
\def\yen{{\mathhexbox@\msafam@55 }}
\def\checkmark{{\mathhexbox@\msafam@58 }}
\def\circledR{{\mathhexbox@\msafam@72 }}
\def\maltese{{\mathhexbox@\msafam@7A }}

\font\tenmsb=msbm10
\font\sevenmsb=msbm7
\font\fivemsb=msbm5
\newfam\msbfam
\textfont\msbfam=\tenmsb
\scriptfont\msbfam=\sevenmsb
\scriptscriptfont\msbfam=\fivemsb
\edef\msbfam@{\hexnumber@\msbfam}

\catcode`\@=\csname pre amssym.def at\endcsname

\expandafter\ifx\csname pre amssym.tex at\endcsname\relax \else \endinput\fi
\expandafter\chardef\csname pre amssym.tex at\endcsname=\the\catcode`\@
\catcode`\@=11
\newsymbol\boxdot 1200
\newsymbol\boxplus 1201
\newsymbol\boxtimes 1202
\newsymbol\square 1003
\newsymbol\blacksquare 1004
\newsymbol\centerdot 1205
\newsymbol\lozenge 1006
\newsymbol\blacklozenge 1007
\newsymbol\circlearrowright 1308
\newsymbol\circlearrowleft 1309
\undefine\rightleftharpoons
\newsymbol\rightleftharpoons 130A
\newsymbol\leftrightharpoons 130B
\newsymbol\boxminus 120C
\newsymbol\Vdash 130D
\newsymbol\Vvdash 130E
\newsymbol\vDash 130F
\newsymbol\twoheadrightarrow 1310
\newsymbol\twoheadleftarrow 1311
\newsymbol\leftleftarrows 1312
\newsymbol\rightrightarrows 1313
\newsymbol\upuparrows 1314
\newsymbol\downdownarrows 1315
\newsymbol\upharpoonright 1316
 
\newsymbol\downharpoonright 1317
\newsymbol\upharpoonleft 1318
\newsymbol\downharpoonleft 1319
\newsymbol\rightarrowtail 131A
\newsymbol\leftarrowtail 131B
\newsymbol\leftrightarrows 131C
\newsymbol\rightleftarrows 131D
\newsymbol\Lsh 131E
\newsymbol\Rsh 131F
\newsymbol\rightsquigarrow 1320
\newsymbol\leftrightsquigarrow 1321
\newsymbol\looparrowleft 1322
\newsymbol\looparrowright 1323
\newsymbol\circeq 1324
\newsymbol\succsim 1325
\newsymbol\gtrsim 1326
\newsymbol\gtrapprox 1327
\newsymbol\multimap 1328
\newsymbol\therefore 1329
\newsymbol\because 132A
\newsymbol\doteqdot 132B
 
\newsymbol\triangleq 132C
\newsymbol\precsim 132D
\newsymbol\lesssim 132E
\newsymbol\lessapprox 132F
\newsymbol\eqslantless 1330
\newsymbol\eqslantgtr 1331
\newsymbol\curlyeqprec 1332
\newsymbol\curlyeqsucc 1333
\newsymbol\preccurlyeq 1334
\newsymbol\leqq 1335
\newsymbol\leqslant 1336
\newsymbol\lessgtr 1337
\newsymbol\backprime 1038
\newsymbol\risingdotseq 133A
\newsymbol\fallingdotseq 133B
\newsymbol\succcurlyeq 133C
\newsymbol\geqq 133D
\newsymbol\geqslant 133E
\newsymbol\gtrless 133F
\newsymbol\sqsubset 1340
\newsymbol\sqsupset 1341
\newsymbol\vartriangleright 1342
\newsymbol\vartriangleleft 1343
\newsymbol\trianglerighteq 1344
\newsymbol\trianglelefteq 1345
\newsymbol\bigstar 1046
\newsymbol\between 1347
\newsymbol\blacktriangledown 1048
\newsymbol\blacktriangleright 1349
\newsymbol\blacktriangleleft 134A
\newsymbol\vartriangle 134D
\newsymbol\blacktriangle 104E
\newsymbol\triangledown 104F
\newsymbol\eqcirc 1350
\newsymbol\lesseqgtr 1351
\newsymbol\gtreqless 1352
\newsymbol\lesseqqgtr 1353
\newsymbol\gtreqqless 1354
\newsymbol\Rrightarrow 1356
\newsymbol\Lleftarrow 1357
\newsymbol\veebar 1259
\newsymbol\barwedge 125A
\newsymbol\doublebarwedge 125B
\undefine\angle
\newsymbol\angle 105C
\newsymbol\measuredangle 105D
\newsymbol\sphericalangle 105E
\newsymbol\varpropto 135F
\newsymbol\smallsmile 1360
\newsymbol\smallfrown 1361
\newsymbol\Subset 1362
\newsymbol\Supset 1363
\newsymbol\Cup 1264
 
\newsymbol\Cap 1265
 
\newsymbol\curlywedge 1266
\newsymbol\curlyvee 1267
\newsymbol\leftthreetimes 1268
\newsymbol\rightthreetimes 1269
\newsymbol\subseteqq 136A
\newsymbol\supseteqq 136B
\newsymbol\bumpeq 136C
\newsymbol\Bumpeq 136D
\newsymbol\lll 136E
 
\newsymbol\ggg 136F
 
\newsymbol\circledS 1073
\newsymbol\pitchfork 1374
\newsymbol\dotplus 1275
\newsymbol\backsim 1376
\newsymbol\backsimeq 1377
\newsymbol\complement 107B
\newsymbol\intercal 127C
\newsymbol\circledcirc 127D
\newsymbol\circledast 127E
\newsymbol\circleddash 127F
\newsymbol\lvertneqq 2300
\newsymbol\gvertneqq 2301
\newsymbol\nleq 2302
\newsymbol\ngeq 2303
\newsymbol\nless 2304
\newsymbol\ngtr 2305
\newsymbol\nprec 2306
\newsymbol\nsucc 2307
\newsymbol\lneqq 2308
\newsymbol\gneqq 2309
\newsymbol\nleqslant 230A
\newsymbol\ngeqslant 230B
\newsymbol\lneq 230C
\newsymbol\gneq 230D
\newsymbol\npreceq 230E
\newsymbol\nsucceq 230F
\newsymbol\precnsim 2310
\newsymbol\succnsim 2311
\newsymbol\lnsim 2312
\newsymbol\gnsim 2313
\newsymbol\nleqq 2314
\newsymbol\ngeqq 2315
\newsymbol\precneqq 2316
\newsymbol\succneqq 2317
\newsymbol\precnapprox 2318
\newsymbol\succnapprox 2319
\newsymbol\lnapprox 231A
\newsymbol\gnapprox 231B
\newsymbol\nsim 231C
\newsymbol\ncong 231D
\newsymbol\diagup 231E
\newsymbol\diagdown 231F
\newsymbol\varsubsetneq 2320
\newsymbol\varsupsetneq 2321
\newsymbol\nsubseteqq 2322
\newsymbol\nsupseteqq 2323
\newsymbol\subsetneqq 2324
\newsymbol\supsetneqq 2325
\newsymbol\varsubsetneqq 2326
\newsymbol\varsupsetneqq 2327
\newsymbol\subsetneq 2328
\newsymbol\supsetneq 2329
\newsymbol\nsubseteq 232A
\newsymbol\nsupseteq 232B
\newsymbol\nparallel 232C
\newsymbol\nmid 232D
\newsymbol\nshortmid 232E
\newsymbol\nshortparallel 232F
\newsymbol\nvdash 2330
\newsymbol\nVdash 2331
\newsymbol\nvDash 2332
\newsymbol\nVDash 2333
\newsymbol\ntrianglerighteq 2334
\newsymbol\ntrianglelefteq 2335
\newsymbol\ntriangleleft 2336
\newsymbol\ntriangleright 2337
\newsymbol\nleftarrow 2338
\newsymbol\nrightarrow 2339
\newsymbol\nLeftarrow 233A
\newsymbol\nRightarrow 233B
\newsymbol\nLeftrightarrow 233C
\newsymbol\nleftrightarrow 233D
\newsymbol\divideontimes 223E
\newsymbol\varnothing 203F
\newsymbol\nexists 2040
\newsymbol\Finv 2060
\newsymbol\Game 2061
\newsymbol\mho 2066
\newsymbol\eth 2067
\newsymbol\eqsim 2368
\newsymbol\beth 2069
\newsymbol\gimel 206A
\newsymbol\daleth 206B
\newsymbol\lessdot 236C
\newsymbol\gtrdot 236D
\newsymbol\ltimes 226E
\newsymbol\rtimes 226F
\newsymbol\shortmid 2370
\newsymbol\shortparallel 2371
\newsymbol\smallsetminus 2272
\newsymbol\thicksim 2373
\newsymbol\thickapprox 2374
\newsymbol\approxeq 2375
\newsymbol\succapprox 2376
\newsymbol\precapprox 2377
\newsymbol\curvearrowleft 2378
\newsymbol\curvearrowright 2379
\newsymbol\digamma 207A
\newsymbol\varkappa 207B
\newsymbol\Bbbk 207C
\newsymbol\hslash 207D
\undefine\hbar
\newsymbol\hbar 207E
\newsymbol\backepsilon 237F
\catcode`\@=\csname pre amssym.tex at\endcsname

\magnification=1200
\hsize=468truept
\vsize=646truept
\voffset=-10pt
\parskip=4pt
\baselineskip=14truept
\count0=1

\dimen100=\hsize

\def\leftill#1#2#3#4{
\medskip
\line{$
\vcenter{
\hsize = #1truept \hrule\hbox{\vrule\hbox to  \hsize{\hss \vbox{\vskip#2truept
\hbox{{\copy100 \the\count105}: #3}\vskip2truept}\hss }
\vrule}\hrule}
\dimen110=\dimen100
\advance\dimen110 by -36truept
\advance\dimen110 by -#1truept
\hss \vcenter{\hsize = \dimen110
\medskip
\noindent { #4\par\medskip}}$}
\advance\count105 by 1
}
\def\rightill#1#2#3#4{
\medskip
\line{
\dimen110=\dimen100
\advance\dimen110 by -36truept
\advance\dimen110 by -#1truept
$\vcenter{\hsize = \dimen110
\medskip
\noindent { #4\par\medskip}}
\hss \vcenter{
\hsize = #1truept \hrule\hbox{\vrule\hbox to  \hsize{\hss \vbox{\vskip#2truept
\hbox{{\copy100 \the\count105}: #3}\vskip2truept}\hss }
\vrule}\hrule}
$}
\advance\count105 by 1
}
\def\midill#1#2#3{\medskip
\line{$\hss
\vcenter{
\hsize = #1truept \hrule\hbox{\vrule\hbox to  \hsize{\hss \vbox{\vskip#2truept
\hbox{{\copy100 \the\count105}: #3}\vskip2truept}\hss }
\vrule}\hrule}
\dimen110=\dimen100
\advance\dimen110 by -36truept
\advance\dimen110 by -#1truept
\hss $}
\advance\count105 by 1
}
\def\insectnum{\copy110\the\count120
\advance\count120 by 1
}

\font\ninerm=cmr9
\font\eightrm=cmr8

\font\tenrm=cmr10 at 10pt

\font\sc=cmcsc10

\def\msb{\fam\msbfam\tenmsb}

\def\bbc{{\msb C}}
\def\bbd{{\msb D}}

\def\bbf{{\msb F}}

\def\bbp{{\msb P}}

\def\bbr{{\msb R}}

\def\bbz{{\msb Z}}

\def\grG{\Gamma}

\def\grO{\Omega}

\def\grS{\Sigma}
\def\grT{\theta}

\def\gra{\alpha}
\def\grb{\beta}
\def\grc{\chi}
\def\grd{\delta}
\def\gre{\epsilon}

\def\grg{\gamma}
\def\gri{\iota}
\def\grk{\kappa}
\def\grl{\lambda}
\def\grm{\mu}
\def\grn{\nu}

\def\grp{\pi}

\def\grs{\sigma}

\def\la#1{\hbox to #1pc{\leftarrowfill}}
\def\ra#1{\hbox to #1pc{\rightarrowfill}}

\def\fract#1#2{\raise4pt\hbox{$ #1 \atop #2 $}}
\def\decdnar#1{\phantom{\hbox{$\scriptstyle{#1}$}}
\left\downarrow\vbox{\vskip15pt\hbox{$\scriptstyle{#1}$}}\right.}

\def\bowtie{\hbox to 1pt{\hss}\raise.66pt\hbox{$\scriptstyle{>}$}
\kern-4.9pt\triangleleft}
\def\hsmash{\triangleright\kern-4.4pt\raise.66pt\hbox{$\scriptstyle{<}$}}
\def\boxit#1{\vbox{\hrule\hbox{\vrule\kern3pt
\vbox{\kern3pt#1\kern3pt}\kern3pt\vrule}\hrule}}

\def\za{\vrule height6pt width4pt depth1pt}

\hfuzz=2pt
\parskip=4pt
\def\intro{1}
\def\ahs{2}
\def\ahpp{3}
\def\stpp{4}
\def\ahpi{5}

\def\tN{\tilde{N}}
\def\bfk{{\bf k}}
\def\bfm{{\bf m}}
\font\svtnrm=cmr17

\def\tX{\tilde{X}}
\def\tN{\tilde{N}}

\def\calo{{\cal O}}

\def\call{{\cal L}}
\def\calm{{\cal M}}

\def\cals{{\cal S}}

\def\Ra*{(R_{a})_{*}}
\def\ada-1{ad_{a^{-1}}~}
\def\p-1u{\pi^{-1}(U)}

\def\PP{{\cal P}{\cal P}}
\def\LPP{{\cal L}{\cal P}{\cal P}}

\def\hol{\hbox{Hol}~(\bbp^1,X)^*}
\def\holk{\hbox{Hol}_{\bfk}(\bbp^1,X)^*}

\centerline{\svtnrm  Stability Theorems for Spaces of Rational Curves}
\medskip

\bigskip
\centerline{\sc C.P. Boyer~~ J.C. Hurtubise~~  R.J. Milgram\footnote{}{\ninerm
During the preparation of this work the first and third authors were
supported by NSF grants and the second author by  NSERC and FCAR grants.}}

\bigskip
\centerline{\vbox{\hsize = 5.85truein
\baselineskip = 12.5truept
\eightrm}}
\centerline{\bf 1. Introduction}
\medskip 

Let  $\grS$ be a compact Riemann surface, and $X$,  a Riemannian manifold.
In the spirit of Morse theory, one would hope that the energy functional
$$E(f) = \int_{\grS}|df|^2 \leqno{\intro.1}$$
on the space of smooth based maps $\hbox{Map}(\grS,X),$ would encode a lot of
the topology of the space in terms of the critical points of the functional.
These critical points are harmonic maps; the absolute minima, in particular,  are often
of  special interest.  It has been known, however, at least  since the work of
Sachs-Uhlenbeck that the Palais-Smale condition C fails for this functional,
making a good Morse theory impossible.  Nevertheless, in some cases, there is a
sense in which the conclusions of Morse theory seem to hold, at least
asymptotically. For simplicity, let us suppose that $\pi_1(X) = 0$, and that
$\pi_2(X)$ is a free abelian group of rank $r$ so that the homotopy class of
any map $f\in \hbox{Map}(\grS,X)$ is given by a multi-degree $\bfk.$ We
can try to compare the homology (homotopy) groups of the space of absolute
minima $\hbox{Min}_{\bfk}(\grS,X)$ with the homology (homotopy) groups of the
full mapping space $\hbox{Map}_{\bfk}(\grS,X).$ If $X$ is K\"ahler, then by a
theorem of Eells and Wood [EW] the space of absolute minima
$\hbox{Min}_{\bfk}(\grS,X)$ is just the space $\hbox{Hol}_{\bfk}(\grS,X)$ of
based holomorphic maps (or anti-holomorphic maps, depending on orientation)
from $\grS$ to $X$, as long as $\hbox{Hol}_{\bfk}(\grS,X)$ is non-empty. We are
then led to consider the inclusion 
$$\hbox{Hol}_{\bfk}(\grS,X) \rightarrow\hbox{Map}_{\bfk}(\grS,X).$$
In
several important cases, one can prove stability theorems, which, loosely
speaking, say that the homology (homotopy) groups of the space of holomorphic
maps is isomorphic to the homology (homotopy) groups of the entire mapping
space through a range that grows with $\bfk$, as $\bfk$ moves to infinity in a
suitable positive cone. These results are compatible with the Morse theoretic
picture, in that in known cases, the indices of the higher critical points also
tend to infinity with $\bfk$, so that if Morse theory were to hold, the
homotopy type of the space of minima would tend to that of the whole space.

Topological stability theorems of the type described had their
origins in questions coming from control theory on the one hand and gauge
theory on the other.  More specifically Atiyah and Jones [AJ] asked whether
such a theorem existed when $\grS$ is replaced by $S^4$ and the functional
instead of \intro.1 is the Yang-Mills functional of instanton gauge theory.
Later it was realized by Atiyah that this was equivalent to  considering a
mapping space $\hbox{Map}(S^2,\call G)$, where $\call G$ is the loop group
corresponding to $G$, the compact Lie group associated to the gauge theory.  A
positive answer to the Atiyah-Jones conjecture was given a few years ago by the
authors and B. Mann in [BHMM2] for the group $G=SU(2),$ following on a 
``stable'' result of Taubes [T] and then more
generally for the various classical compact Lie groups by Tian [Ti1, Ti2] and
Kirwan [Ki2].  The conjecture should hold for more general four-manifolds 
than $S^4$; see [HM] for the case of a ruled complex surface. However the
first proof of any such stability theorem was given in the mapping space
(``sigma model'')  case by Segal [Se] for $X=\bbp^n$.  Segal's theorem was then
generalized to the case $X$ a complex Grassmannian by Kirwan [Ki1], and certain
$SL(n,\bbc)$ flags by Guest [Gu1].  Later, Mann and Milgram [MM1] increased the
range of the isomorphisms obtained by Kirwan for Grassmannians and treated
[MM2] all $SL(n,\bbc)$ flag manifolds.  Moreover, the essential technique of
L-stratifications used by the authors in the present and previous papers
[BHMM1,BHMM2] was introduced in [MM1,MM2].  More recently, 
topological stability theorems were proven for any generalized flag manifold
$G/P$ by the authors [BHMM1], [Hu1], following on a stable result of 
Gravesen [Gra] and  for toric varieties by Guest
[Gu2]. We refer to [Hu2] for a survey.

In all of the above cases an essential ingredient is
representing the minima $\hbox{Min}_{\bfk}(\grS,X)$ by a labelled configuration
space.  This description is essentially confined to spaces which admit actions
of groups of dimension equal to that of the space. In what follows we prove the
stability theorem for what is more or less the largest class of manifolds to
which this particle description applies. It includes, for example, all the
preceeding (non-singular) cases, as well as  many smooth compact  spherical varieties
(for definitions, see section 2 below). It is
interesting to note that, with the exception of the Atiyah-Jones conjecture for
which the target space $X$ is infinite-dimensional, all other varieties for
which the theorem was known to hold were spherical; the present paper covers
many $X$ for which this is not the case, simple examples being provided by blow-ups 
of $\bbp^n$ along varieties in the hyperplane at infinity. One is still left, however, with the
following important question:

\noindent{\sc Question}: What is the most general complex target space $X$
which admits a topological stability theorem of the type described above?

For such a theorem to hold, the space of based holomorphic maps should 
build up in a regular way the entire space of based maps.  The space
$\hbox{Map}(\grS,X)$ is often quite complicated topologically; thus,
 the spaces $\hbox{Hol}_{\bfk}(\grS,X)$ of
parameterized based holomorphic curves of a particular genus on $X$, must have,
at least stably, an equivalent structure.  In particular,  for  $\grS = \bbp^1
= S^2$ the Riemann sphere, if the loop space $\Omega^2(X)$ is non trivial,
 there should be a plethora of rational curves on $X$, passing through, for example, any 
finite set of points. Such varieties are called  rationally connected [Ko].
The role of rational curves in varieties is appearing to be more and more essential
in their study, whether 
it is to define quantum cohomology, or to study higher dimensional varieties in Mori's 
program (in which the rationally connected varieties play an 
important role) (See e.g. [Mi,Ko]). Our theorem can be seen as another
manifestation of this, in a more homotopy theoretic setting; 
 it is interesting that the topological property of stability should appear to 
be linked to some sort of rationality property of the manifold.

We thus 
 prove the stability results for 
a large  subclass of rationally connected varieties, in fact rational varieties,
on which a complex solvable linear group acts with  a free dense open orbit.   We
call such varieties {\it principal almost solv-varieties}. Actually, 
 we can do a
little better, and remove the condition that the group be linear, 
 by noticing that stability holds trivially for Abelian varieties,
and making use of the Albanese map $\gra: X\ra{1.3} \hbox{Alb}(X).$  It turns
out that in our case $\gra$ is a locally trivial fibration whose fibres $F$ are
principal almost solv-varieties, and that the long exact sequence in homotopy implies
that $\pi_2(X)\simeq \pi_2(F).$  We will  show that $\pi_2(X)$ is free,
and that it embeds into a fixed lattice $\bbz^\kappa.$ The multi-degree of an
element of a homotopy class is then given in terms of this embedding. Our first
main result is:

\noindent{\sc Theorem} A: [Homology  Stability Theorem] \tensl Let $X$ be a
smooth compact K\"ahler manifold which has a holomorphic action of a connected
complex solvable Lie group $S$ with an open orbit $N$ on which $S$ acts freely.
Then for every multi-degree $\bfk = (k_1,\cdots,k_\kappa)$ the inclusion $\gri({\bf
k}(X))$ induces an isomorphism in homology with $\bbz$ coefficients through
dimension $q({\bf k})$ between the space of holomorphic maps of degree $\bfk$
and the space of continuous maps of degree ${\bf k}$:  
$$(\gri({\bf k}(X)))_t: H_t(\hbox{Hol}_{\bf k}(\bbp^1,X)^*;\bbz) ~\cong~
H_t(\grO^2_{\bf k}X;\bbz)$$
for 
$$t~\leq~q({\bf k})~=~ c_0 l({\bf k}) -1.$$ 
Here  $c_0(X)>0 $ is 
a constant  depending only on the space $X$,  and
$l({\bf k})= min~(k_i)$.  \tenrm

 To prove this theorem, we define stabilisation
maps from  the 
spaces of holomorphic maps of degree ${\bf k}$ to  spaces of  maps of higher degree
${\bf k'}$.
The proof of  theorem A uses a result of Gravesen
[Gra], which gives a homology isomorphism between  the loop space
and the limit of the 
spaces of holomorphic maps. Theorem A tells us that the 
stabilisation of the homology occurs in the nicest possible way, so that one does not, 
for instance, have for a fixed homology group, an infinite sequence of classes 
being created in  moduli spaces of maps of ever higher degrees, only to 
disappear as one increases the degree a bit 
further.

In particular cases, the constant  $c_0(X) $ can be estimated; some comments will be 
given after the proof of  theorem A in section 4 below.
 The controlling factor is the structure of the singular 
set of the complement $X_\infty$ of the free dense orbit of the solvable group. If,
for example, $X_\infty$ has smooth components with normal crossings (e.g.
$\bbp^n$, toric varieties, Bott-Samelson varieties), one  can obtain
$c_0(X)  = 1$; for  cases such as $\bbp^n$,  this is of course 
still far from the   
$c_0(X) =2n-1$ obtained  by Segal [S]; it is however 
the known range for toric varieties [Gu2].

Since $\hbox{Hol}_{\bf k}(\bbp^1,X)^*$ and $\grO^2 X$ are not always
simply connected, our second main result, Theorem B below, does not trivially
follow from Theorem A. Indeed, in general, the gap between homology 
isomorphisms and homotopy isomorphisms can be quite large, for example
if one  considers homology spheres versus homotopy spheres. 

\noindent{\sc Theorem} B: [Homotopy  Stability Theorem]  \tensl
Let $X$ be as in Theorem A. 
Then for all {\bf k} the inclusion
$\gri({\bf k}(X))$ is a homotopy equivalence through dimension
$c_0 l({\bf k}) -s -2,$ where $c_0$ and $l(\bfk)$ are as in Theorem A, and
$s$ is the rank of $\pi_3(X)= \pi_1(\grO^2 X).$  \tenrm

The proof of theorem B, which is regrettably but necessarily quite
technical since it must give information on the structure of the homology
of cyclic covers of quite complicated spaces, the strata at infinity
described above, is given in Sections 5 and 6.

The authors would like to thank Paulo Lima-Filho for information
on spherical varieties, Bruce Gilligan for communicating Akhiezer's
examples to us, and Ben Mann for helpful discussions.

\bigskip

\noindent{\S2. \bf Principal Almost Solv-varieties} 
\def\rkX{\hbox{rk}~(X)}

\medskip

\noindent{\it (i) Almost homogeneous spaces} 

Let $(X,\calo_X)$ be an irreducible complex space, and let $G$ be a connected
complex Lie group acting as biholomorphic maps of $(X,\calo_X).$  
 The definition given below is due to Remmert and Van de Ven [cf.
Akh,Huc,HO].

\noindent{\sc Definition} \ahs.1: \tensl Let $X$ be an irreducible complex
space and let $G$ be a complex Lie group acting as biholomorphic
transformations on $X.$  We say that $X$ is {\it almost homogeneous} if $G$
has an open orbit $N$ in $X.$  If $G$ acts freely on $N$ so that $N$ can be
identified with $G$ itself, we call $X$ a {\it principal almost homogeneous
space}.  \tenrm

Notice that a given complex space $X$ can be almost homogeneous with respect to
different Lie groups or with respect to different actions of the same Lie
group. When possible confusion can arise we will specify the Lie group $G$
and/or the open orbit $N.$  

The complement $X\backslash N$ of $N$ in $X$ is an
analytic subspace of $X$ [Akh], and  $X$ can be viewed as an equivariant
compactification of $N$ by adding $X\backslash N$.  We shall refer to
$X\backslash N$ as ``infinity'' in $X,$ and denote it by $X_{\infty}.$ (When $X=\bbp^n$,
with the action of $\bbc^n$, $X_\infty$ is indeed the hyperplane at infinity; when one switches
to the group $(\bbc^*)^n$, $X_\infty$ is the union of the coordinate hyperplanes). The Lie
group $G$ acts biholomorphically on the irreducible components of $X_{\infty}$.
As a variety $X_{\infty}$ may have many irreducible components.  For example
for $X=G/P$ a generalised flag manifold, with the action of the unipotent group
``opposite" to $P$,  one finds that $X_{\infty}$ is the union of the closures of
all codimension one Bruhat cells in the Bruhat decomposition of $G/P.$  The
number of irreducible components of $X_{\infty}$ then equals the number of
codimension one Bruhat cells. When $P$ is a Borel subgroup, this is the rank $r$
of the group; more generally, one has $s\le r$ components,  with $s$ equal to
 the dimension of the center of the
Levi factor of the parabolic subgroup $P.$

Our first result provides a wealth of examples of almost homogeneous spaces.
The proof of the following proposition is straightforward.

\noindent{\sc Proposition} \ahs.2: \tensl Let $X$ be an almost homogeneous
space with respect to the Lie group $G,$ and let $V$ be a complex subvariety
of $X_{\infty}$ of codimension greater than 1 in $X$ that is stable under the
action of $G.$  Let $\pi:\tilde{X}\ra{1.4} X$ denote the blow-up of $X$ along
$V.$ Then $\tX$ is almost homogeneous with respect to $G$, with open orbit
$\tN$ given by the inverse image of $N.$ In particular, if $N$ is a principal
orbit in $X$ so is $\tN$ in $\tX.$ \tenrm

We shall be especially interested in certain types of almost homogeneous
spaces, namely those where the Lie group in question is a solvable complex Lie
group $S.$ 

\noindent{\sc Definition} \ahs.3: \tensl An irreducible complex space $X$ that
is almost homogeneous with respect to a connected solvable complex Lie group
$S$ is called an {\it almost solv-space}.  Furthermore, if the action of $S$ is
free on the open orbit, $X$ is called a {\it principal almost solv-space}. If
$X$ is a complex manifold, we use the terminology {\it almost solv-manifold}.
If $X$ is an algebraic variety and $S$ is a connected solvable linear algebraic group
acting algebraically on $X$, then $X$ is called an {\it almost solv-variety}.
\tenrm

Examples are given below.

\noindent{\sc Proposition} \ahs.4: \tensl Let $X$ be a normal almost
solv-variety.  Then infinity $X_{\infty}$ has pure codimension one in $X;$
hence, $X_{\infty}$ is a Weil divisor on $X.$  \tenrm

\noindent{\sc Proof}: By a theorem of Snow [Huc] the open orbit $N\simeq
\bbc^k\times (\bbc^*)^l$ as a complex manifold. Hence, $N$ is holomorphically
convex.  Thus, for each component of $X_{\infty}$, there are holomorphic
functions on $N$ which do not extend holomorphically across it.  If this
component had codimension $\geq 2,$ Hartog's theorem for normal varieties
[Iit,pg 124] would imply that every holomorphic function would extend across
it.  Thus, every component of $X_{\infty}$ has codimension $1$ in $X.$
\hfill\za

\noindent {\sc Proposition \ahs.5}\tensl\ Any principal almost solv-variety
is rational.\tenrm

\noindent{\sc Proof:} As above let $n = k+l$ be the complex dimension of $X$.
Since $N\simeq
\bbc^k\times (\bbc^*)^l$, we have embeddings as dense open sets
$X\hookleftarrow N\hookrightarrow \bbp^n $
which defines a birational map. \hfill\za

\noindent{\it (ii) Spherical varieties}

\noindent{\sc Definition} \ahs.6: \tensl Let $G$ be a connected complex
reductive algebraic group.  A normal algebraic variety $X$ on which $G$ acts as
a group of algebraic transformations is called a {\it spherical variety} if some
Borel subgroup $B$ of $G$ has a dense orbit in $X.$  Let $H$ be the stabiliser in $G$
of a point in the dense orbit; the   homogeneous space $G/H$   is an open dense
subvariety of $X$. The subgroup
$H\subset G$ is called {\it spherical}.   \tenrm

  It
is clear that a spherical variety is an almost solv-variety with respect to the
connected solvable linear complex group $B.$ One can, however, often find a smaller solvable 
subgroup which has a {\it free} dense orbit.
To see this we need to use the structure of spherical homogeneous
spaces due to Brion, Luna, and Vust [BLV]. Let $G,H,B$ be as in Definition
\ahs.6. Then $BH$ is open in $G.$ Let $P$ denote the subgroup of $G$ defined by
$P=\{p\in G ~|~ pBH=BH\}.$  Since $P$ contains $B$ it is parabolic, and thus
has a Levi decomposition $P\simeq L\bowtie U$ where $L$ is reductive and
$U=R_u(P)$ is the unipotent radical.  Let $C=R(L)$ be the radical of $L.$ It is
an algebraic torus and the connected component of the center of $L.$  Let $N$
denote the dense open $B$-orbit in $G/H\subset X.$ A theorem of [BLV] (see also [Br]) says that
as varieties $N\simeq U\times (C/C\cap H).$ Consider the exact sequence of
algebraic groups 
$$<e>\ra{1.5}C\cap H\ra{1.5} C\ra{1.5} C/C\cap H\ra{1.5} <e>.$$ 
If $C\cap H$ is connected then this sequence splits [Bo1]. Let
$T\subset C$, $T\simeq C/C\cap H$ denote the image of $C/C\cap H$ under such a
splitting.  For any such choice of $T$ we denote by $S$ the subgroup of $B$
which is the semi-direct product $T\bowtie U.$ Then by construction $S$ acts
freely on $N\simeq U\times (C/C\cap H).$ We have arrived at:

\noindent{\sc Proposition} \ahs.7: \tensl Let $G$ be a reductive algebraic
group, and let $H,B,C,T$ be subgroups of $G$ as described above. Suppose
further that $C\cap H$ is connected.  Then $S =T\bowtie U$ acts freely on the
dense open $B$-orbit $N.$ Hence, every spherical variety $X$ with $C\cap H$
connected is a principal almost solv-variety with respect to the connected
solvable linear algebraic group $S.$ \tenrm

\noindent{\sc Remark} \ahs.8: The following simple examples due to Akhiezer shows that
the connectivity of $C\cap H$ in Proposition \ahs.7 is indeed necessary. First
take $G=SL(2,\bbc)$ and $H=\bbc^*$ any torus in $G.$ Then any solvable subgroup
$S$ of $G$ with a dense orbit must be a Borel subgroup for  dimensional reasons. Moreover,
the intersection of $S$ with any torus contains the center $=\bbz_2,$ so $S$
never acts freely. In this case $S=B$ and $C$ is its reductive factor, and one
easily sees that $C\cap H =\bbz_2.$ For  another example let $G$ be any
complex semi-simple Lie group, and take $H$ to be an extension of a maximal
unipotent subgroup by a finite group. Then the intersection $C\cap H$ will be a
non-trivial finite group. 

For spherical varieties, there is a notion of rank, defined (see Brion [Br]) to be the
dimension of $C/C\cap H$ in the proof of Proposition \ahs.7.  
More generally,  we have

\noindent{\sc Definition} \ahs.9: \tensl We define the {\it rank} of any
principal almost solv-variety $X,$ denoted by $\hbox{rk}~(X),$ to be the
dimension of the torus $S/S_u,$ where $S$ is the connected solvable algebraic
group of Definition \ahs.3 and $S_u$ is the subgroup of unipotent elements of
$S.$ \tenrm

The rank is in a sense the crudest invariant of a principal almost
solv-variety. Even in the spherical case a given spherical variety can be
given the structure of a principal almost solv-variety in many inequivalent
ways. This will be   illustrated  below.
 
\noindent{\it (iii) Examples of principal almost solv-varieties}

\noindent {\sc Generalized flag manifolds} \ahs.10:  Let $G$ be a semi-simple
complex Lie group, and let $P$ be a parabolic subgroup of $G$.  The quotient
$G/P$ is a generalised flag manifold, and is obviously spherical.  The
unipotent radical $U$ of $P$ acts freely on an open orbit of $G/P$, as do all
the conjugates of $U$. One can choose this conjugate to be ``opposite'' to $P$,
so that its Lie algebra is a sum of those root spaces not lying in the Lie
algebra of $P$. A variety $X$ is a complete spherical variety of $\rkX =0$ if
and only if the spherical subgroup is parabolic and $X=G/P$ [Br]. 

\noindent {\sc Toric varieties} \ahs.11: By definition a toric variety is a
variety together with an action of $(\bbc^*)^n$ which is free on an open dense
orbit.  These varieties are spherical, with $G= (\bbc^*)^n$, and $H = <e>,$ and
principal almost solv-varieties. Moreover, $X$ is a toric variety if and only
if $\rkX = \hbox{dim}~X.$

\noindent{\sc $\bbp^n$ as a principal almost solv-variety} \ahs.12: 
Some varieties can be principal solv-varieties in a variety of distinct ways.
This is the case for projective spaces $\bbp^n$. Indeed, one can
obtain them as flag manifolds (rank 0), toric varieties (rank $n$), or indeed,
any rank between zero and $n$. For $\bbp^2$, one has for example a whole
infinite sequence of rank one solvable groups acting in affine coordinates
by
$$T_{\lambda, \zeta} (x,y) = (\lambda^k x, \lambda^{l}y + \zeta),$$
where $k\neq 0$ and $k,l$ are coprime or $l=0$. These are, in fact, 
all the ways to obtain $\bbp^2$ as a rank 1 almost solv-variety, up to isomorphism.

\noindent {\sc Some examples of complete symmetric varieties} \ahs.13: Let $\sigma$ be the
 involution of 
$PGL(n,\bbc)$ induced by conjugation by the matrix
$$\pmatrix {0&I_r&0\cr I_r&0&0\cr 0&0&I_s},$$
$n=2r+s$, and let $H$ be the fixed point subgroup of $\sigma$. The variety 
$SL(n,\bbc)/H$
is spherical, and De Concini and Procesi [DP] construct a spherical compactification $X$ 
of  this space which is smooth. One can show that the subgroup $P$  is the parabolic subgroup 
whose Lie algebra is the sum of root spaces with roots $\alpha$ satisfying
$$<(r,r-1, ...,1,-r,-r+1,...,-1,0, 0,...,0), \alpha>\ \ge 0,$$ and that $C\cap H$ is connected.

Next we have examples of principal almost solv-varieties which are not
spherical:

\noindent{\sc Equivariant Blow-ups} \ahs.14: As noted in \ahs.2, 
one can take any $S$-invariant subvariety of codimension 
greater than one in a principal almost-solv variety $X$ and blow it up to obtain
another principal almost-solv variety $\tilde X$.  The invariant subvariety must then lie 
in $X_\infty$.  Taking $X$ to be $\bbp^n$, for example, one can blow up any subvariety of the
plane at infinity and obtain a principal  almost-solv variety for the group $\bbc^n$.
We note that this blowing up can reduce the automorphism groupp; for example, blowing up 
 up a set of 
points in $\bbp^n$ reduces
the automorphism group from $PGL(n,\bbc)$ to the subgroup which fixes that set.
Blowing up a suitably generic set of points along infinity in $\bbp^n$ thus reduces the 
automorphism group to the $\bbc^n$ of translations, which acts trivially 
along infinity. Such a variety is definitely not spherical. More generally, 
one obtains a whole set of examples from the standard ones by blowing up subvarieties which are
invariant under the action of the solvable group. 

\noindent {\sc Bott-Samelson varieties} \ahs.15: Let us consider a given
reductive group $G$, a fixed maximal torus $T,$ and a Borel subgroup $B,$ with
negative simple roots $s_1,...,s_r$ in the Lie algebra of $G$ corresponding to
generators $w_1, w_2,...,w_r$ of the Weyl group.  For each $i$, we have the
parabolic subgroups $P_i$ whose Lie algebra is obtained by adjoining the  root
space of $s_i$ to the Lie algebra of $B$, so that $P_i/B= \bbp^1$. For each
word $\tilde w= w_{i_1}w_{i_2}\cdots w_{i_k}$ in the generators of the Weyl
group, we consider the quotient 
$$P(\tilde w)= P_{i_1}\times P_{i_2}\times \cdots \times P_{i_k}/ B^k,$$
where $B^k$ acts by $(b_{i_1}, b_{i_2},\dots,b_{i_k})
(p_{i_1},p_{i_2},\dots,p_{i_k}) =
(p_{i_1}b_{i_1}^{-1},b_{i_1}p_{i_2}b_{i_2}^{-1},\dots,b_{i_{k-1}}p_{i_k}
b_{i_{k}}^{-1})$. This is smooth, and when the word is a reduced expression for
an element $w$ for the Weyl group, it is a desingularisation of the Schubert
variety $S_w$ given as the closure of the orbit $BwB/B$ in $G/B$. Over the open
cell, the map to $BwB/B$  is an isomorphism.  The cell $BwB/B$ is a free orbit
of a unipotent subgroup $U_w$ of $B$, and this action lifts to $P({\tilde w})$,
giving it the structure of a principal almost-solv variety.  One can show that
the divisor at infinity has smooth components with normal crossings.  These
varieties are frequent objects of study in representation theory [Jan].

\noindent {\sc Products} \ahs.16 One can of course, take products of almost-solv-varieties
to obtain new ones. For example, one obtains the  compact
complex symmetric spaces (in the sense of Borel [Bo2]) associated to a reductive
group $G$  by taking  products of flag 
manifolds and certain toric varieties ([Leh])

\noindent{\it (iv) The topology of principal almost solv-varieties}

For any rational variety $X,$ in fact for unirational varieties, there is a
vanishing theorem $H^q(X,\calo)=0$ for $q>0.$  In addition, if $X$ is
smooth,  the exponential exact sequence gives for the line bundles on $X$:
$$\hbox{Pic}(X)\simeq H^2(X,\bbz). \leqno{\ahs.17}$$
Next, returning to our decomposition
$X=N\cup X_{\infty},$ we write $X_{\infty}$ as the union of irreducible
algebraic hypersurfaces $X_{\infty}=\cup_{\gra =1}^\grk X_{\gra},$ and have: 

\noindent{\sc Proposition} \ahs.18: \tensl Let $X$ be a connected smooth
principal almost solv-variety, of complex dimension $n$. Then $X$ is simply connected, $H_2(X,\bbz)$ is
torsion free and $\hbox{Pic}(X)\simeq H^2(X,\bbz)\simeq H_2(X,\bbz)\simeq
\bbz^{\grk-l}$ where $\grk$ is the number of irreducible components of
$X_{\infty}$ and $l= \hbox{rk}~(X).$ \tenrm

\noindent{\sc Proof}: By Proposition \ahs.5 a principal almost solv-variety
is rational, and the fundamental group is a birational invariant for algebraic
varieties [GH, pg 494]. Thus, $X$ is simply connected.

The action of the solvable linear algebraic group $S$ on the pair $(X,N)$
gives the commuting diagram, with exact rows, in homology with integer coefficients: 
$$\matrix{\matrix{\cdots &H_*(S)\otimes H_*(N)&\fract{id\times
\gra_*}{\ra{1.5}}& H_*(S)\otimes H_*(X)&\fract{id\times j}{\ra{1.5}}&
H_*(S)\otimes H_*(X, N)\cr
&\decdnar{\grb_*}&&\decdnar{\grb_*}&&\decdnar{\grb_*}\cr
\cdots&H_*(N)&\fract{\gra _*}{\ra{1.5}}&H_*(X)&\fract{j}{\ra{1.5}}&H_*(X,N)\cr}
\cr \hbox to 1in{\hfill}\matrix{ \fract{id\times \partial}{\ra{1.5}}&
H_*(S)\otimes H_{*-1}(N) &\fract{id\times\gra_{*-1}}{\ra{1.5}}& \cdots\cr
&\decdnar{\grb_*}\cr
\fract{\partial}{\ra{1.5}}&H_{*-1}(N)&\fract{\gra_{*-1}}{\ra{1.5}}&
\cdots\cr}\hfill\cr} \leqno{\ahs.19}$$ 
where $\grb_*$ is the induced map in homology.  Under our assumptions, 
$H_*(N)= H_*(S)$, and $H_*(S)$ is  generated by $H_1(S)$.

In dimension one the simple connectivity of $X$ implies that
$\partial \colon H_2(X, N)\ra{1} H_1(N) = \bbz^l$ is onto.  
Since a set of generators for $H_i(N)$ has the form $\{ \grb_*(e_j\otimes
v_j)\}$ where the $e_j$ run over $H_{i-1}(S)$ and the $v_j$ generate $H_1(N)$,
the formula
$$\partial \grb_*(e_j\otimes w_j) = -\grb_*(e_j\otimes\partial(w_j))
\leqno{\ahs.20}$$
shows that $\partial\colon H_i(X, N) \ra{1} H_{i-1}(N)$ is onto for each $i$.
This implies that in each dimension greater than or equal to two we have the
short exact sequence
$$0\ra{1.5} H_i(X)\ra{1.5} H_i(X,N)\ra{1.5} \bbz^{l_{i-1}}\ra{1.5} 0
\leqno{\ahs.21}$$
where $H_{i-1}(N) = \bbz^{l_{i-1}}$, and consequently, splittings
$$H_i(X, N) = H_i(X) \oplus \bbz^{l_{i-1}}.$$
By Alexander-Poincare
duality 
$$H_2(X,N)\simeq H^{2n-2}(X_{\infty})\simeq \bbz^\grk.$$ 
So the exact sequence \ahs.21 shows that $H_2(X)\simeq \bbz^{\grk-l}.$ In
particular $H_2(X)$ is torsion free, and so by the universal coefficients theorem,
$H_2(X)\simeq H^2(X).$ \hfill\za

Let $r$ denote the rank of $\hbox{Pic}(X)\simeq H_2(X)\simeq \pi_2(X).$ The
components of $\grO^2X$ are labeled by a multi-degree $\bfk =(k_1,\cdots,k_r)$
and $\grO^2_{\bfk}X$ denotes such a component.  There is a basis
$\{L_i\}_{i=1}^r$ for $\hbox{Pic}(X)$ such that for $f\in \grO^2_{\bfk}X,$
$c_1(f^*L_i)=k_i$ for all $i=1,\cdots,r.$  It is also convenient to give a
description in terms of the line bundles $[X_{\gra}]$ associated to the divisor
$X_{\gra}$ with certain relations.  Dualizing the exact sequence \ahs.21 gives
$$\matrix{0&\ra{1.5}&H^1(N)&\ra{1.5}& H^2(X,N)&\ra{1.5}& H^2(X)&\ra{1.5}&0\cr
&&\decdnar{\simeq}&&\decdnar{\simeq}&&\decdnar{\simeq}&& \cr
0&\ra{1.5}&\bbz^l&\ra{1.5}&H_{2n-2}(X_{\infty})&\fract{\psi}{\ra{1.5}}
&\hbox{Pic}(X)&\ra{1.5}&0.  \cr} \leqno{\ahs.22}$$
The map $\psi$ sends the class $\hbox{cl}(X_{\gra})$ in
$H_{2n-2}(X_{\infty})$ representing the irreducible component $X_{\gra}$
to the line bundle $[X_{\gra}].$  As this is an epimorphism of free Abelian
groups, for each $L_i$ there is a sum $V_i$ of components $X_{\gra}$ such that
$\psi(\hbox{cl}(V_i))=L_i.$ Thus, $L_i$ and $[V_i]$ define the
same element of $\hbox{Pic}(X),$  so $c_1(f^*[V])=k_i.$ Alternatively,
$k_i$ is the intersection number of $f(S^2)$ with the variety $V_i$ at
infinity. 

\bigskip

\noindent{\S3. \bf The Poles and Principal Parts Description}

\medskip

\noindent{\it (i) The sheaf of principal parts}

For any complex space $X$ we let ${\cal O}(X)$ denote the
sheaf of germs of holomorphic maps from $\bbp^1$ into $X$, and
${\cal O}(U,X)$ denote the holomorphic sections of ${\cal O}(X)$
over the open set $U\subset \bbp^1$. If the space $X$ has a base
point, ${\cal O}(X)$ will be a sheaf of pointed sets.  We denote the global
holomorphic maps by $\hbox{Hol}~(\bbp^1,X),$ and the global based holomorphic
maps by $\hbox{Hol}~(\bbp^1,X)^*$. There is natural inclusion
$$\hbox{Hol}~(\bbp^1,X)^*\hookrightarrow \grO^2 X \leqno{\ahpp.1}$$
into the space of all based continuous maps $\grO^2X.$    

Let $X$ be an almost solv-variety, with $N=X\backslash X_\infty$ the dense
open orbit in $X.$ Following Gravesen [Gra] we define the presheaf of
meromorphic maps from $\bbp^1$ to $N$ by setting for each open $U$ in $\bbp^1$:
$${\cal M}(U,N)= {\cal O}(U,X)\backslash{\cal O}(U,X_{\infty}) \leqno{\ahpp.2}$$
and we let ${\cal M}(N)$ denote its associated sheaf. Let $S$ be the connected
solvable linear algebraic group that acts biholomorphically on $X$ and
transitively on $N.$ There is a natural action of the sheaf $\calo(S)$ on
$\calm(N)$ given on the presheaf level by pointwise multiplication: for each
open set $U\subset \bbp^1,$ the action  of $\calo(U,S)$ on $\calm(U,N)$  sends
a local section $f_U\in \calm(U,N)$ to the local section $g_U\cdot f_U$ defined
by $g_U\cdot f_U(z) = g_U(z)\cdot f_U(z).$ This allows us to define the
quotient sheaf 
$$\PP \simeq \calm(N)/\calo(S) \leqno{\ahpp.3}$$
to be the sheaf of germs of equivalence classes of
meromorphic maps where $f,f'\in \calm(U,N)$ are equivalent if there is $g\in
\calo(U,S)$ such that $f'=g\cdot f.$ The sheaf $\PP$ is called the sheaf of
{\it principal parts} in $X.$  Notice
that the stalk $\PP_z$ is actually independent of the location $z\in \bbp^1$ of the pole:
one can simply translate the maps from point to point, by an automorphism 
of $\bbp_1$; over $\bbc\subset \bbp^1$, this can be done unambiguaously by translations.
We thus  denote the space $\PP_z$ of ``local principal parts'' simply by $\LPP$.

To give an idea of what the stalk of this sheaf 
looks like, we give a few examples. We consider the  stalk at $z=0$.

1) When $X=\bbp^n, N= \bbc^n$ and $S=\bbc^n$, the principal parts 
construction yields exactly the classical principal parts of a meromorphic
map. Indeed, one has that the local form of a map into $\bbp^n$ near $z=0$ 
is given by an n-tuple of Laurent series. Let us write the $i$-th entry
of this $n$-tuple as $a_kz^k + a_{k+1}z^{k+1}+...$. Normalising
these series under the additive action of maps into $\bbc^n$ just allows us to
kill of all the positive order terms, leaving just the negative order terms,
i.e., the principal part. We note that in our sense a principal part encodes
both some discrete information (the order of the pole) as well as some
continuous information (the principal part or Laurent tail) 

2) When $X= \bbp^n, N= (\bbc^*)^n$ and $S = (\bbc^*)^n$,
one obtains the other classical description of a map in terms of zeroes in
homogeneous coordinates. Indeed,
one can normalise the series $a_kz^k + a_{k+1}z^{k+1}+...$ under the multiplicative action of
maps into $(\bbc^*)$ to the monomial $z^k$, so that the information contained in the 
principal part is  simply the order of each component, or, equivalently, the multiplicity 
of intersection with each component of $X_\infty$, which in this case are the coordinate
hyperplanes. The information contained in the principal part is then  discrete
in nature.

3)  The same calculation of a normal form holds for the maps into a
 general toric variety: indeed, one has the same local coordinates as in 2) above.
Thus the information contained in a principal part in the toric case
is simply the multiplicity of intersection with the components
of the  divisor at infinity, which in this
case are simply the closures of the codimension one orbits. This ties in with the description
of the maps given in [Gu2].

4) A computation of some spaces of principal parts for the space of full flags 
in $\bbc^3$ is given 
in [BHMM1].

The sheaf $\PP$ in many ways resembles the sheaf of divisors on a complex manifold, in that
when one considers a fixed section and restricts it to stalks, one obtains 
the trivial element (corresponding to maps into $N$) generically, 
and, over exceptional points (the {\it support} of the section), 
some non-trivial element.  We will, in all cases, refer to the points in the support
as {\it poles}, even though in some examples, such as $\bbp_1$ in its toric description,
some of the poles are actually zeroes! 

\noindent{\it (ii) Configurations of principal parts}

Now assume that $X$ is a smooth almost solv-variety.  An element $P\in
H^0(\bbp^1,\PP)$ can be represented by a sequence of pairs $(U_i,f_i)$ where
$\{U_i,i=0,\cdots,n\}$ is a finite cover of $\bbp^1$ and $f_i\in \calm(U_i,N).$
The pull-back $f_i^*X_{\infty}$ is a divisor on $U_i,$ i.e. a finite sum of
points.  Furthermore, we can choose the cover $\{U_i\}$ so that $U_i\cap U_j =
\phi$ for $i,j>1$, $f_0(U_0)\subset N$, and $f_i(U_i),i=1,\cdots,n$ intersects
infinity $X_{\infty}$ at a single point $f_i(z_i).$ We shall refer to such
covers as $P$-{\it good} covers. A global section in $H^0(\bbp^1,\PP)$ is
called a {\it configuration of principal parts}. It consists of a finite number
of points $z_i\in \bbp^1$ (location of the poles), together with the {\it
local} principal parts data, a non-trivial element in the stalk
$\LPP =\PP_{z_i}\simeq {\cal M}_{z_i}/{\cal O}(S)_{z_i}$ at each point $z_i$. 
We view $H^0(\bbp^1,\PP)$ as  a labelled configuration space with $\LPP$ as the
space of labels. The topology on $H^0(\bbp^1,\PP)$ is the quotient topology
from $\calm(U,N)$ with the later given the compact-open topology.  See [Gra]
for details.

An element in $\hbox{Hol}~(\bbp^1,X)^*$ will naturally determine
a section in $H^0(\bbp^1,\PP)$, whose poles are the points mapped to
$X_{\infty}$. We now fix the basing condition.   We choose the base point in
$\bbp^1\simeq S^2\simeq \bbc\cup \{\infty\}$ to be the north pole $\{\infty\}$
and the base point of $X$ to be a fixed point $*\in N\subset X$. This precludes
a map in $\hbox{Hol}~(\bbp^1,X)^*$ having a pole at $\infty\in \bbp^1$.  Let
$H^0(\bbc,\PP)$  denote the subspace of $H^0(\bbp^1,\PP)$ whose poles are all
located in $\bbc\simeq \bbp^1 -\{\infty \}$, so that our holomorphic map
determines an element of $H^0(\bbc,\PP)$.  There is then a map
$$\hol\ra{1.5} H^0(\bbc,\PP), \leqno{\ahpp.4}$$ 
which we shall see is an inclusion, so that
the holomorphic maps are determined by their principal parts. 
We are interested in when a configuration of
principal parts $P\in H^0(\bbc,\PP)$ comes from a holomorphic map $f\in \hol.$
In general there is an obstruction cocycle coming from the  action of
$(\bbc^*)^l\subset S$ which we shall describe shortly. 

\noindent{\it (iii) Virtual multi-degrees and label spaces}

We have just seen that a based holomorphic map $f\in \hbox{Hol}~(\bbp^1,G/P)^*$ corresponds to a
finite set of points $f^*X_{\infty}\subset \bbp^1$ together with labels in a space
$\LPP.$ These labels carry information about the map. One piece of information 
the labels always carry is   that of a virtual multi-degree. 

Indeed, any element $LP$ of $\LPP$ is represented by a map $f$ of a disk in $\bbc$ into 
$X$, with the center of the disk being mapped to $X_\infty$ and the rest of the 
disk being sent to $N= X\backslash X_\infty$. The  maps of the disk
into $S$ acts on this by homotopies, and so one has a well-defined multiplicity
 or {\it virtual degree} of $LP$ as an element of
$H_2(X,N;\bbz)$ which from the proof of Proposition
\ahs.16 is $\bbz^\grk.$ There exist natural generators for $H_2(X,N;\bbz)$,
given by the Alexander Poincar\'e duals to the
components $X_\alpha$ of the divisor $X_\infty$. In this basis,
we write the multiplicity  as a vector
$$ (m_1,\dots,m_{\grk})\leqno \ahpp.5$$
Analytically, we can obtain the $m_\gra$ by considering the holomorphic sections 
$s_{\gra}$ of the line bundles
$[X_{\gra}]$ vanishing at the hypersurface $X_{\gra}.$ These are   invariant
up to a factor by the action of $S$.  The multiplicity $m_\gra$
is the order of the vanishing of the pull-back $f^*s_\gra$. 

 Let
$P$ belong to $ H^0(\bbc,\PP)$, with poles $z_i,i=1,\cdots,r$, and local
representatives $f_i:U_i\rightarrow X$ for the principal parts. We have multiplicities
associated to each $f_i$, and so, summing, a multiplicity or virtual degree for $P$:

$$\bfk(P) =\sum_{i=1}^r{\bf m}^i =(k_1,\cdots,k_{\grk}) \qquad  k_{\gra}=
\sum_{i=1}^{r}m^i_{\gra}.  \leqno{\ahpp.6}$$ 

We  define some notation. If the point $z_i$ is mapped by $f_i$ to the divisor $X_\alpha$,
we will refer to $z_i$ as an $\alpha$-pole. We note that $z_i$ can be 
both an $\gra$-pole and an $\grb$-pole, by being mapped to $X_\gra\cap X_\grb$.
 An $\gra$-pole is called {\it
simple} if $m_\grb =1~\hbox{for $\grb =\gra$}$ and $m_\grb =0$ if $\grb\neq
\gra.$ The label $\alpha $ will sometimes be referred to as the {\it color} of the pole.
We can also
define the {\it scalar multiplicity} of the $i^{th}$ pole by $$|{\bf m}^i|=
\sum_{\gra =1}^{\grk}m^i_{\gra}. \leqno{\ahpp.7}$$ Similarily, the scalar
virtual degree is $|{\bf k}|= \sum_\gra k_\gra.$

We can decompose  $H^0(\bbc,\PP)$  as the disjoint union
$$H^0(\bbc,\PP)=\bigsqcup H_{\bf k}^0(\bbc,\PP),\leqno{\ahpp.8}$$ where $H_{\bf
k}^0(\bbc,\PP)$ is the subspace of those principal parts with virtual
multi-degree equal to $\bfk.$ In the same way, the local principal parts space
decomposes according to multi-degree:
$$\LPP = \bigsqcup_{\bf m} \LPP_{\bf m},\leqno \ahpp.9$$ 

\noindent{\sc Example} \ahpp.10: To fix our ideas, we will consider a simple example,
that of degree one based maps from $\bbp^1$ to
$\bbp^2 $, to see what configurations we can obtain.  The based holomorphic map
 $f:\bbp^1\ra{1.3} \bbp^2$
are given in homogeneous coordinates by $f(z)=[z-z_0,z-z_1,z-z_2].$ The only
constraint for this to be a well defined map is that all homogeneous
coordinates cannot vanish simultaneously. Thus the space of degree one maps
in $\hbox{Hol}~(\bbp^1,\bbp^2)^*$ is $\bbc^3-\{z_0=z_1=z_2\}.$ The poles and 
principal parts description of this space varies according to the choice of group.

\item{(1)} When one considers $\bbp^2$ with the action of translation of $\bbc^2$,
 infinity is a single hyperplane
divisor, say $X_0$, defined by setting the first homogeneous coordinate equal to
zero.  There is one pole, with position $z_0$. Fixing the pole, the
label space $\LPP_1$ corresponds to the residues $z_0-z_1$ and $z_0-z_2,$ and so
$\LPP_1\simeq \bbc^2-\{(0,0)\}$. The space of maps is then an $\LPP_1$ fibration over 
$\bbc$.  

\item{(2)} At the opposite extreme is the toric description of $\bbp^2.$ Here
infinity is the union $\cup_{\gra=0}^2X_{\gra}$ of the coordinate hyperplane
divisors $X_{\gra}=\{[x_0,x_1,x_2] |x_{\gra}=0 \},\gra= 0,1,2$. There
are three ``poles'' of different colors, with position $z_0,z_1,z_2$  with the
constraint that all three $z_i$ are not equal, corresponding to the fact that the 
three coordinate lines do not intersect. The spaces of degree one maps splits
into strata: the generic stratum corresponds to the $z_i$ being distinct,
so that there is one pole of multiplicity $(1,0,0)$, one of multiplicity
$(0,1,0)$, and one of multiplicity $(0,0,1)$. The other 3 strata correspond to
two of the poles coinciding, so that there is one stratum of maps with one pole
of multiplicity $(1,1,0)$ and one of multiplicity $(0,0,1)$, one stratum with
multiplicities $(1,0,1), (0,1,0)$ and one stratum with multiplicities $(0,1,1),
(1,0,0).$ In all cases,  the label spaces $\LPP_{(i,j,k)}$ are single points.

\item{(3)} The intermediate case corresponds to  the rank one examples of 
\ahs.12. In this case infinity $X_{\infty}$ is the
singular variety $\{[x_0,x_1,x_2]|x_0x_1=0\}.$ There are two poles $z_0, z_1$,
and the mapping space divides into two strata, according to whether $z_0, z_1$
coincide or not.  When the poles are different, one has one pole of multiplicity
$(1,0)$, with label in $\LPP_{(1,0)}=pt$ and one pole of multiplicity
$(0,1)$, with label in $\LPP_{(0,1)}=\bbc$.
When the poles coincide, one then has one pole of multiplicity $(1,1)$, and the 
label space $\LPP_{(1,1)}$ is $\bbc^*$.

The ``additivity'' of principal parts plays a crucial role in stability
theorems. Let $P,P'\in H^0(\bbc,\PP)$ be such that they have no poles in
common, then their union $P\cup P'$ is a configuration of principal parts in an
obvious way. Let us take an isotopy of $\bbc$ into a disk $D$; this defines a 
diffeomorphism $H^0(\bbc,\PP) \simeq H^0(D,\PP)$, (which is not holomorphic).
This allows us to define on $H^0(\bbc,\PP)$ the structure of a
homotopy-associative monoid: one first maps the configurations of
principal parts over $\bbc$ into configurations of principal parts
over disjoint disks $D$ and $D'$, then adds them.  This gives 

\noindent{\sc Proposition} \ahpp.11: \tensl There is a continuous map
$$A:H_{\bf k}^0(\bbc,\PP)\times H^0_{\bf l}(\bbc,\PP)\ra{1.5} 
H^0_{{\bf k}+{\bf l}}(\bbc,\PP), $$
where ${\bf k}+{\bf l}$ has components $k_{\gra}+l_{\gra}$. \tenrm

\noindent

\noindent{\it (iv) Configurations and holomorphic maps}

It is quite straightforward to see that the map $\hol\ra{1.3} H^0(\bbc,\PP)$
is injective. Indeed, one has that any two based maps with the same
principal parts differ by the action of a global map $\bbp^1\rightarrow S$; 
these maps are constant, and the basing condition guarantees that the constant
is the identity, so that the maps coincide.

We  now consider the  question of which configurations of principal parts
correspond to holomorphic maps. 
When a configuration of  principal parts $P$ is obtained from  a
based holomorphic map $f\in \hbox{Hol}~(\bbp^1,X)^*,$ the corresponding
multi-degree $\bfk$   is constrained to lie in the image of $H_2(X;\bbz)
\subset H_2(X,N;\bbz)$.  We will show that this is the only 
obstruction to a configuration of principal parts being derived
from a map.   

As with principal parts we denote the subspace of
$\hbox{Hol}~(\bbp^1,X)^*$ of maps having multi-degree $\bfk$ by
$\hbox{Hol}_{\bfk}(\bbp^1,X)^*.$   For any complex space $M$ with a sheaf of
groups $\calo(S)$ we recall the definition of $H^1(M,\calo(S)).$ Given a cover
$\{U_i\}$ of $M$ by open sets we consider holomorphic maps $g_{ij}:U_i\cap
U_j\ra{1.3} \calo(S)$ which satisfy the cocycle condition $g_{ij}g_{jk}g_{ki}=
e$ in $U_i\cap U_j\cap U_k.$ Two such maps $g_{ij}$ and $g'_{ij}$ are
equivalent if for each $i$ there is a holomorphic map $\phi_i:U_i\ra{1.3}
\calo(S)$ such that for all $i,j$ we have $g'_{ij}=\phi_i^{-1}g_{ij}\phi_j.$
Then $H^1(M,\calo(S))$ is defined to be the limit, under refinement of the cover,
of the  set of equivalence classes of such
equivalence maps and  it  classifies equivalence classes
of $S$-bundles over $M$ . If the cover is acyclic (Leray), one can compute $H^1(M,\calo(S))$
directly from the cover, without taking limits, as in the Abelian case. 

Let us now consider our case, where $X$ is a principal almost homogeneous space
with respect to the complex Lie group $S.$ We recall that $S$
acts freely on the open orbit $N,$ so that $N$ can be identified with $S$
itself.  Given a configuration of principal parts $P$ represented by local
holomorphic maps $f_i\in \calm(U_i,N)$ with respect to the $P$-good cover
$\{U_i\},$ we consider the restrictions of $f_i$ and $f_j$ to the intersection
$U_i\cap U_j.$ Since $S$ acts freely on $N,$ there is a unique holomorphic map
$g_{ij}:U_i\cap U_j\ra{1.3} \calo(S)$ such that $f_i=g_{ij}\cdot f_j$ in
$U_i\cap U_j.$ Moreover, by refining the cover if necessary we can check that a
different representative $f'_i\in \calm(U_i,N)$ of $P$ gives rise to a cocycle
$g'_{ij}$ that is related to $g_{ij}$ by $g'_{ij}= g_ig_{ij}g^{-1}_j$ where
$g_i:U_i\ra{1.3} \calo(S)$ are holomorphic.  So we get a well-defined element
$[g_{ij}]\in H^1(\bbp^1,\calo(S)).$ Principal parts arising from a holomorphic map
give the trivial cocycle, and so  we have a sequence of maps of
pointed sets
$$<e>\ra{1.5} \hbox{Hol}~(\bbp^1,X)^*\ra{1.5}
H^0(\bbc,\PP)\fract{\grd}{\ra{1.5}} H^1(\bbp^1,{\cal O}(S)), \leqno{\ahpp.12}$$
and we would like to think of this as being exact. 

\noindent{\sc Definition} \ahpp.13 [Gro]: \tensl Let $A,B,C$ be
sets with $C$ a set with a neutral element $e.$ Then the sequence
$$A\ra{1.5} B\fract{\nu}{\ra{1.5}} C$$ 
is said to be {\it exact} (at $B$) if the kernel $\nu^{-1}(e)$
equals $\hbox{image}~(A).$ \tenrm

\noindent Warning: If $A=0$ (a one point set), then
exactness at $B$ does not imply the injectivity of $\nu.$  We refer the reader to the lecture
notes of Grothendieck [Gro] for further detail. 

In our case,  the set
$H^1(\bbp^1,{\cal O}(S))$ has a preferred or ``neutral'' element, namely the
class of the identity element in $C^1(U_i\cap U_j,\calo(S)).$ 
 To obtain exactness, note that an element of 
$H^0(\bbc,\PP)$ that maps to a trivial cocycle in $H^1(\bbp^1,{\cal O}(S))$
can be split   as
$\phi_i\phi_j^{-1}.$  As this then guarantees that the local maps $f_i$ patch
together to give a global section in $H^0(\bbp^1,\calm(N))^*$, we have exactness, that is 
 the following Mittag-Leffler type result:

\noindent{\sc Theorem} \ahpp.14: \tensl  Let $X$ be a smooth complex principal
almost homogeneous space.  Then the configurations of principal parts (i.e.
elements of $H^0(\bbc,\PP)$) that can be represented by based holomorphic maps
(i.e elements of $\hbox{Hol}~(\bbp^1,X)^*$) are precisely those lying in the
kernel of the map $\grd$ of \ahpp.12.  \tenrm

Let us return to the case at hand when $S$ is a connected solvable algebraic
group.  There is a short exact sequence of sheaves of groups on $\bbp^1,$
$$<e>\ra{1.5}\calo(R)\ra{1.5}\calo(S)\ra{1.5}\calo(T)\ra{1.5} <e>,
\leqno{\ahpp.15}$$ 
where $T$ is an algebraic torus of rank $l$ and $R$ is the normal subgroup of
unipotent elements of $S.$ Since $R$ is unipotent we have a vanishing
theorem for $H^1(\bbp^1,\calo(R)).$ When $R$ is Abelian and unipotent this is
quite classical, since $\calo(R)$ is isomorphic to a sum of copies of the
structure sheaf $\calo.$ However, when $R$ is not Abelian this is not so well
known, so we give the proof; see also [Gra].

\noindent{\sc Lemma} \ahpp.16: \tensl $H^1(\bbp^1,\calo(R))= 0.$ \tenrm

\noindent{\sc Proof}:  We have the composition series 
$<e>=R_n<R_{n-1}<\cdots <R_0=R$ for some positive integer $n,$ with each factor
$R_i/R_{i+1}$ Abelian. So we have the exact sequence
$$<e>\ra{1.5} \calo(R_{i+1})\ra{1.5} \calo(R_{i})\ra{1.5}
\calo(R_{i}/R_{i+1})\ra{1.5} <e>$$ 
of sheaves of groups on $\bbp^1.$ According to section 5.3 of Grothendieck
[Gro] this gives an exact sequence (as in Definition \ahpp.14)
$$H^1(\bbp^1,\calo(R_{i+1}))\ra{1.5} H^1(\bbp^1,\calo(R_{i}))\ra{1.5}
H^1(\bbp^1,\calo(R_{i}/R_{i+1})).$$ 
But since $R_n=<e>,$ $R_{n-1}$ is Abelian, and each factor group $R_i/R_{i+1}$
is Abelian, the first and last terms of this sequence vanish.  But this implies
the vanishing of the middle term.  An easy induction finishes the proof.
\hfill\za

We are now ready for:

\noindent{\sc Proposition} \ahpp.17: \tensl Let  $X$ be a principal almost
solv-variety. Then there is a map
$H^1(\bbp^1,\calo(S))\fract{\nu}{\ra{1.3}} H^1(\bbp^1,\calo(T))\simeq
\bbz^l.$ Furthermore, $\gamma\in H^1(\bbp^1,\calo(S))$ maps to zero in
$\bbz^l\simeq H^1(\bbp^1, {\cal O}(T))$ if and only if it represents the
trivial bundle.  \tenrm

\noindent{\sc Proof}: The exact sequence \ahpp.15 implies the existence of the
map $\nu$ (See [Gro], section 5.3). If $\grg$ represents the trivial bundle in 
$H^1(\bbp^1,\calo(S)),$ one easily sees that $\nu(\grg)=0.$ Conversely, suppose
that $\grg$ maps to $0$ in $H^1(\bbp^1,\calo(T)),$ and choose a
good cover $U_i$, so that the only non empty double intersections are $U_0\cap
U_i, i=1,\cdots,n$. Let $s_{0i} \in C^1(U_0\cap U_i,\calo(S))$ be a cocycle
representing $\gamma$. We can use the fact that (\ahpp.15) splits to write
$s_{0i} = u_{0i}t_{0i}$ with $u_{0i}\in C^1(U_0\cap U_i,\calo(R))$, and
$t_{0i} \in C^1(U_0\cap U_i,\calo(T))$.  By hypothesis, one can write $t_{0i}=
t_0 ^{-1}t_i$, so that $s_{0i}$ can be replaced by the equivalent cocycle
$\tilde s_{0i} = t_0s_{0i}t_i^{-1} = t_0(u_{0i}t_0^{-1}t_i)t_i^{-1} =
t_0u_{0i}t_0^{-1}$. Since the subgroup $R$ is normal, this last cocycle lies
in $\calo(R)$, and therefore, by the vanishing theorem \ahpp.16, can be
written as $u_0u_i$, expressing $\tilde s_{0i}$, and therefore $s_{0i}$, as a
coboundary. \hfill\za

We recall that a principal part gives a well defined element of the relative 
homology $H_2(X,N)$. Also, there is a well defined map from $H^1(\bbp^1, {\cal
O}(S))$ to $H_1(N;\bbz)$, given with respect to any good cover of  $\bbp^1$ by 
$U_0 = \bbp^1-\{p_1,...,p_j\}$ and disjoint  disks $U_i $ centred at $p_i$, as follows:
  A cocycle $g_{0i}$ for this covering can be
restricted to   circles around the punctures, giving a well defined element of
$H_1(N;\bbz)$, under the identification of $N$ with $S$. The map 
$H^1(\bbp^1, {\cal
O}(S))\rightarrow H_1(N;\bbz)$ factors through $H^1(\bbp^1, {\cal O}(T))$, and in
fact gives an isomorphism $H^1(\bbp^1, {\cal O}(T) \simeq H_1(N;\bbz))$.
With this
identification, we have an exact sequence of commutative diagrams:
$$\matrix {<e>&\ra{1.5} &\hbox{Hol}~(\bbp^1,X)^*&\ra{1.5}
&H^0(\bbc,\PP)&\fract{\grd}{\ra{1.5}}& H^1(\bbp^1,{\cal O}(S))&\ra{1.5} &0 \cr
& &\downarrow\phi & & \downarrow\psi& & \downarrow\rho&\cr
0&\ra{1.5}&H_2(X) &\ra{1.5} & H_2(X,N) &\ra{1.5} &
H_1(N)&\ra{1.5}&0}.\leqno{\ahpp.18} $$

The map
$\psi $ associates to each configuration of principal parts a 
virtual multi-degree $(k_1,\cdots, k_\kappa)$ in $H_2(X,N)$. The 
 ``real'' or ``good'' multi-degrees $(k_1,\cdots, k_\kappa)$ are those 
in the subgroup $H_2(X)\simeq\bbz^{\kappa-l}\subset
H_2(X,N)$.  From the diagram, one has:

\noindent{\sc Proposition} \ahpp.19: \tensl A configuration of principal parts
$P\in H^0(\bbp^1,\PP)$ represents 
 a holomorphic map $f\in \hol$
if and only if the virtual multi-degree ${\bf k}=(k_1,\cdots,k_\kappa)$ maps to
zero in $H_1(N).$ For multi-degrees which do map to zero, one has a
homeomorphism $$H^0_{\bfk}(\bbc,\PP)\simeq \holk.$$ \tenrm

In particular, this gives the principal part spaces the structure of
algebraic varieties.

\vfil\eject

\bigskip

\noindent{\S\stpp. \bf Homology Stabilization of Principal Parts}

\medskip

\noindent{\it (i) A stratification of $H^0_\bfk(\bbc,\PP)$}

In this section unless otherwise stated $X$ will denote a smooth 
 principal almost solv-variety with respect to a connected solvable
linear algebraic group $S.$  We will exploit the
isomorphism given by proposition \ahpp.19 to give a proof of homology stability
theorems. The proof relies on two key phenomena.
The first is that spaces of labelled particles, such as configurations of principal
parts, exhibit a natural stability property in homology as one increases the number of 
particles.
 Thus one must check  that the
space behaves sufficiently like a space of many particles. This is ensured by
the second key point, proving that the space of maps with poles of high multiplicity (and therefore
fewer poles) has a high homological codimension. This in turn is done by
showing that the space is smooth, and that the geometric codimension of the
high multiplicity set is high.

\noindent{\sc Proposition} \stpp.1: \tensl Let $X$ be a smooth complex
principal almost solv-variety, and $\bfk$ a good virtual multi-degree.
Then the space $H^0_{\bfk}(\bbc,\PP)=\holk$ is a smooth finite dimensional
complex manifold.  \tenrm

\noindent{\sc Proof}: The deformation theory proof of smoothness of $\hol$
in Theorem 4.1 of [BHMM1] carries over without any essential changes to our
case.  This is due to the fact that the transitivity of the action of $G$
on the dense affine open set $N$ guarantees that at a map $f$ the pull-back
$f_i^*TX$ splits as the sum of positive line bundles, i.e,
$$f_i^*TX\simeq \oplus_i({\cal O}(\lambda_i)) \leqno{\stpp.2}$$
with $\grl_i \geq 0.$ The obstruction to smoothness $H^1(\bbp^1,f^*TX) $ then
vanishes. \hfill\za

Let us fix an irreducible component $X_{\gra}$ of $X_\infty.$  We let $H^0_{{\bf
k},\alpha}(\bbc,\PP)$ denote the subspace of $H^0_{\bf k}(\bbc,\PP)$ of those
principal parts with only one $\gra$-pole, which is then  of multiplicity
$k_{\gra}$; the other poles can be of arbitrary multiplicity.  More generally, for any $S$-stable subvariety $V\subset X_{\gra}$,
we define $H^0_{{\bf k},\alpha}(\bbc,\PP)(X,V)$ to be the subset of $H^0_{{\bf
k},\alpha}(\bbc,\PP)$ consisting of those principal parts $P$ whose $\gra$-pole
have their image in $V$.  In particular, we are interested in the cases when
$V=X_{\gra}^{*}$, the smooth locus at infinity, or
$V=X_{\gra}^s=X_{\gra}-X_{\gra}^*$, the singular locus at infinity. Also, if
${\bf k}$ is a multi-index, we let $H^0_{{\bf k,k}}(\bbc,\PP)$ be the
subvariety of $H^0_{{\bf k}}(\bbc,\PP)$ with only one pole of multiplicity
${\bf k}$, i.e.  $\gra$-multiplicity $k_{\gra}$ for each $\gra
=1,\cdots \grk$.  We remark that $H^0_{{\bf k,k}}(\bbc,\PP)$ could be empty.

\noindent{\sc Proposition} \stpp.3:  \tensl Taking ``codimension'' to mean 
complex codimension in $H^0_{{\bf k}}(\bbc,\PP)$, 
\item{(i)}The codimension of $H^0_{{\bf k},\alpha}
(\bbc,\PP)(X,X^{*}_{\gra})$ is $k_{\gra}-1$.
\item{(ii)}  The codimension of $H^0_{{\bf
k},\alpha}(\bbc,\PP)(X,X^{s}_{\gra})$ is bounded below by
$max~\{C_{\gra}k_{\gra} -C_{\gra}^{\prime},1\}$. 
\item{(iii)} There are constants
$0<C_{\gra}\leq 1$ and $C^{\prime}_{\gra}\geq 0$ depending only on $\gra$ and 
$X$ such that the codimension of $H^0_{{\bf k},\alpha}(\bbc,\PP)$ 
 is bounded below by $C_{\gra}k_{\gra} -C^{\prime}_{\gra}$, and by $1$ if $k_\alpha>1$.
\item{(iv)} There are constants
$0<\hat C\leq 1$ and $\hat C^{\prime}\geq 0$ depending only on $X$ such that
the codimension of $H^0_{{\bf k,k}}(\bbc,\PP)$ is bounded below by $\hat
C|{\bf k}| -\hat C^{\prime}$. 
\item{(v)} Unless $k_\alpha = \delta_{\alpha, \beta}$ for some $\beta$,
the codimension of  $H^0_{{\bf k,k}}(\bbc,\PP)$ is at least one. \tenrm

\noindent{\sc Proof}:   The
proof  differs only slightly from the proof of Theorem 5.20
of [BHMM1], given there for the generalized flag manifolds $G/P$. We summarise here
the salient points.
We begin with the proof of (i). Let us place ourselves at a point $P$ of
$H^0_{{\bf k},\alpha}(\bbc,\PP)(X,X^{*}_{\gra})$, with $\alpha$-pole at $z$. We stabilise 
$H^0_{{\bf k},\alpha}(\bbc,\PP)(X,X^{*}_{\gra})$ by adding in
a large number of  principal parts 
  located away from the support of $P$, in such a way that 
the final multidegree ${\bf k+l}$ corresponds to a space of holomorphic maps. We then have 
a map $H^0_{{\bf k}}(\bbc,\PP)
\rightarrow H^0_{{\bf k+l}}(\bbc,\PP)$, mapping $P$ to a configuration $P'$
corresponding to a map $f'$. One has a local 
product structure $H^0_{{\bf k+l}}(\bbc,\PP)\simeq H^0_{{\bf k}}(\bbc,\PP)\times
H^0_{{\bf l}}(\bbc,\PP)$. The codimension of $H^0_{{\bf k},\alpha}(\bbc,\PP)(X,X^{*}_{\gra})$
at $P$ is the same as the codimension at $P'$ of the subvariety of 
$H^0_{{\bf k+l}}(\bbc,\PP)(X,X^{*}_{\gra})$ consisting of maps with one $\alpha$-pole 
of multiplicity $k_\alpha$ and all its other $\alpha$-poles simple. 

One then can use the fact that the pull-back $f'^*(TX)$ is sufficiently positive
to show that the evaluation map  taking elements of $H^0_{{\bf k+l}}(\bbc,\PP)$
to their $k_\alpha$-jets at $z$, is submersive at $P'$. Again, details are given in 
5.20 of [BHMM1]. 

The problem is reduced to understanding tangencies in the space of jets.
One shows that the jets tangent to $X_\alpha$ to 
order $k_\alpha$ have codimension $k_\alpha$; this tells us that
the maps 
having at least one point with a  jet tangent to $X_\alpha$ to 
order $k_\alpha$ have codimension $k_\alpha-1$. As we are on the smooth locus, 
we can choose local coordinates $x_1,...x_n$, with $X_\alpha$ cut out by $x_1 = 0$.
The desired result is then a simple observation.

For (ii), one must also analyse jets at singular points of the hypersurface $X_\alpha$.
 It then
appears that the codimension of the subspace of jets vanishing to order $k_\alpha$
 is not necessarily  $k_\alpha$, but can be lower. By blowing up, however, 
one can obtain the estimate 
$C_\alpha k_\alpha
-C'_\alpha$. Again, we refer to [BHMM1]. A general position argument to note that the generic
holomorphic map does not meet the singular locus of $X_\alpha$, 
showing that the codimension is at least one. Part (iii) follows from (i) and (ii).

Parts (iv) and (v) are obtained by repeating the arguments of (i), but with the divisor
$X_\infty$ instead of its components $X_\alpha$.
\hfill\za

Referring to \ahpp.9, we can let  $\LPP_{\bf k}$ be the subspace of $H^0_{{\bf k,k}}(\bbc,\PP)$ with the pole
at $0$, so that $H^0_{{\bf k,k}}(\bbc,\PP) = \LPP_{\bf k}\times \bbc$.

As in [BHMM1] we consider a collection of multi-indices
$${\cal M}~=~ \{ {\bf m}^1, \dots, {\bf m}^r \} \leqno{\stpp.4} $$
satisfying the constraint \ahpp.6. 

\noindent{\sc Definition} \stpp.5: \tensl Let ${\cal V}_{{\cal M}}$  be
the subset of all elements of $H^0_{{\bf k}}(\bbc,\PP)$ that have
$r$ poles at   distinct points $z_1,\dots,z_r$ of $\bbc$ such that the
multiplicity of the pole at $z_i$ is ${\bf m}^i$.  \tenrm

Now we have the ``pole location map'' 
$$\Pi: H^0_{{\bf k}}(\bbc,\PP)\ra{1.5} SP^{|\bfk|}(\bbc). \leqno{\stpp.6}$$
defined by sending a principal part $P$ to the unordered set of points of
$\bbc$ defining the configuration of poles counting multiplicity. That is, if
$z_i\in \bbc$ is a pole of multiplicity $\bfm^i$ then $z_i$ occurs $|\bfm^i|$
times.The map $\Pi$ is not necessarily
surjective, as for example, there can be a sequence of divisors at infinity with empty
intersection. Thus, not all poles can necessarily coalesce. Nevertheless, the
pole location map \stpp.6 restricted to ${\cal V}_{\calm}$ gives the locally trivial
fibration
$$\Pi_{{\cal M}}:{\cal V}_{{\cal M}} \ra{1.5} \bbd\bbp^{r}(\bbc) \leqno{\stpp.7}$$ 
where $\bbd\bbp^r(\bbc)$ denotes the deleted $r$-fold product; that is, the
space of $r$ distinct unordered points in $\bbc$. The fiber of $\Pi$ in \stpp.7
is $\prod_{i=1}^r \LPP_{{\bf m}^i}$ which is not necessarily smooth. Thus, the
strata ${\cal V}_{\cal M}$ are not necessarily smooth; however, they are varieties, as
they are cut out of $H^0_{{\bf k}}(\bbc,\PP)$ by the vanishing of derivatives
of locally defined analytic functions. Furthermore, one can subdivide the
spaces $\LPP_{{\bf m}^i}$ into smooth strata; this induces a subdivision
${\cal V}_{{\cal M},i}$ of the ${\cal V}_{\cal M}$. One can order the ${\cal M}$ in such a
way that
 together with Proposition \stpp.3
one has

\noindent{\sc Theorem} \stpp.8:  \tensl There is a stratification of
$H^0_{\bfk}(\bbc,\PP)$ by smooth complex varieties ${\cal V}_{\calm,j}$ which
satisfies: 
\item{(i)} $\dim~{\cal V}_{{\cal M},j} <\dim~{\cal V}_{{\cal M}',i}$
only if ${\cal M}' < {\cal M}$, or, if ${\cal M}'={\cal M}$, then $i< j$,
\item{(ii)} ${\cal V}_{{\cal M},j}\subset \overline{{\cal V}_{{\cal M}',i}}$
only if
 ${\cal M}' \leq {\cal M}$, or,
if ${\cal M}'={\cal M}$, then $i\leq j$
\item{(iii)} The set $\{\calm,j\}$ is finite.

\noindent Furthermore, there is a positive constant $c(X)$, which is
independent of the stratum indices ${\cal M},j$ and the multi-degree ${\bf k}$,
so that the complex codimension of ${\cal V}_{{\cal M},j}$ in $H^0_{{\bf
k}}(\bbc,\PP)$, $$\hbox{codim}({\cal M},j;{\bf k}) ~=~ \hbox{codim}({\cal V}_{{\cal
M},j}\subset H^0_{{\bf k}}(\bbc,\PP)) $$ is bounded below by $$c(X)\sum_i
(|{\bf m}^i|-1) ~=~ c(X)(|{\bf k}|-r)$$ where $r$ is the number of poles.
Here, referring to proposition \stpp.3,  $c(X)$ is 
a constant such that  $c(X)n-1 \le {\rm max}(1,\hat C n-\hat C')$ for all $n$.
 \tenrm

This last fact follows from \stpp.3 and the local decomposition of 
$H^0_{{\bf k}}(\bbc,\PP)$ near ${\cal V}_{\cal M}$ as a product $\Pi_{i=1}^r
H^0_{{\bf m^i}}(\bbc,\PP)$. 
Thus the big strata are the ones with many poles of low multiplicity. As in
[BHMM1] we construct a stabilization map as follows: For each irreducible
component $X_{\gra}$ fix a local principal part of multiplicity one. We have 

\noindent {\sc Proposition} \stpp.9 \tensl There is a rational map of
multi-degree $(k_1,\cdots,k_\kappa)$ with all $k_i>0$. \tenrm

\noindent {\sc Proof:} By Proposition \ahs.5 the variety $X$ is birational to
$\bbp^n$, and maps $f:\bbp^1\rightarrow X$ get transformed into maps $\tilde
f:\bbp^1\rightarrow \bbp^n$.  The constraint that, say, $0\in \bbp^1$ gets
mapped into $X_\alpha$ gets transformed into a constraint on $\tilde f(0),
\tilde f'(0), \tilde f''(0),\cdots,\tilde f^{(s)}(0)$, for some $s$. However,
maps into $\bbp^n$ are given by polynomials, and one can arrange by
interpolation for a map to satisfy any finite set of constraints such as those
given above, at a set of points $z_1, z_2,\cdots, z_k$ of $\bbp^1$. In other
words, one can ensure that the map $f$ meets each $X_\alpha.$ \hfill\za

\noindent{\it (ii) Stabilization of principal parts}

Proposition \stpp.9 tells us that the plane of good virtual degrees
$(k_1,\cdots
k_\kappa)$ intersects the interior of the positive quadrant
$k_i>0$.  Let us choose some small good  multi-degree ${\bf l} $ in this
quadrant, and build a configuration $P_{\bf l}$ of multi-degree ${\bf l}$ with
simple poles from our local principal parts.  Using Proposition \ahpp.8, we add
$P_{\bf l}$ to an arbitrary $P\in H^0_{\bfk}(\bbc,\PP).$ This gives a
continuous map $$\gri({\bf k},\bfk +{\bf l}): H^0_{\bf k}(\bbc,\PP) \ra{2}
H^0_{\bfk + {\bf l}}(\bbc,\PP)\leqno{\stpp.10}$$ This map preserves the
stratification, so that it can be studied stratum by stratum.  The arguments in
[BHMM1] then carry over word for word to give:

\noindent{\sc Theorem} \stpp.11: \tensl For all ${\bf k}$
 the inclusion \stpp.10
induces an isomorphism in homology 
$$(\gri({\bf k},\bfk +{\bf l}))_t: H_t(H^0_{\bf k}(\bbc,\PP);\bbz)
~\cong~ H_t(H^0_{\bfk +{\bf l}}(\bbc,\PP);\bbz)$$
for 
$$t\leq q~=~q({\bf k})~=~
c_0(X)l({\bf k}) - 1,$$
where $c_0(X)>0$ 
is a constant  which depends only on $X$ and where
$l({\bf k})= \min(k_1,\dots,k_{\kappa}).$  \tenrm

We recall that the proof proceeded by considering the stabilisation map
stratum by stratum; the homology of the whole is then obtained from 
a Leray type spectral sequence, whose $E^1$-term is  the homology of the strata
suspended by the real codimension of the stratum. Ordering the strata as in Theorem
\stpp.8, labelling them as ${\cal V}_{{\cal M}_1} \ge {\cal V}_{{\cal M}_2}\ge...$,
one has $E^1_{ij}= H_{i+j-{\rm codim}} ({\cal V}_{{\cal M}_i})$
 The fact that the stabilisation maps preserves
the spectral sequence gives a map between the $E^1$ terms:
$$E^1( H^0_{\bf k}(\bbc,\PP))_{ij}\longrightarrow E^1(H^0_{\bfk +{\bf l}}(\bbc,\PP))_{ij}$$
which commutes with differentials. The isomorphism on the homology of the whole
spaces for a given range follows from the isomorphism on the constituent pieces,
through the same range, plus one.

On the stratum ${\cal V}_{{\cal M}_i}$, the homology stabilises in a range which is 
half the minimum (over the set of colors) of the numbers
of simple poles in ${\bf m}$. (The map in homology is injective on the strata; 
one then simply has to look  for the smallest class one can build in the charge ${\bf m+l}$
stratum which one cannot build in the charge ${\bf m}$ stratum.) To the stability range of the homology 
map on the stratum one must then add the real codimension of the 
stratum, which is bounded from below, using \stpp.8, by $2(c(X)(|{\bf k}|-r)  )$.
The range for the stability of the whole space is then $c_0(X)l({\bf k}) - 1$,
with $c_0 = {\rm min}(2c(X), 1/2)$.

Using the maps $\gri(\bfk,\bfk + {\bf l})$ of \stpp.10 one can form the
direct limit $\displaystyle{\lim_{\ra{1}}H^0_{\bf k}(\bbc,\PP)}$ and obtain
natural inclusions
$$\gri({\bf k}): H^0_{\bf k}(\bbc,\PP) \ra{1.5} \lim_{\ra{1}}H^0_{\bf
k}(\bbc,\PP). \leqno{\stpp.12}$$  

Theorem \stpp.11 tells us then that the  map $\gri({\bf k})$ induces isomorphisms
in homology in the range $t\le c_0|{\bfk}| - 1$. We now use the
argument of Gravesen [Gra], showing how the limit space is homologically 
equivalent to $\Omega^2(X)$. This argument, given there for maps to 
the $SL(n,\bbc)$ flag manifolds only,  depends on having
a pole and principal parts description, so it holds equally well in
the more general case of a principal almost solv-variety. Indeed, the argument 
of Gravesen, in essence, works inductively,  building up on one side the 
continuous maps
of the 2-sphere (or, more generally, any surface)
 into $X$ from maps of disks into $X$, and using natural fibration 
properties of spaces of maps under restriction. On the particle space side, the 
same procedure is used, building up the space of particles over the sphere from 
particles over disks, with quasi-fibrations replacing fibrations in 
the restriction maps. At the last step in the induction (adding the last disk
 and closing off the surface), one has to  take the
mapping telescope $\widetilde{H^0(}\bbc,\PP)$ of the stabilized principal part
space $\displaystyle{\lim_{\ra{1}}H^0_{\bf k}(\bbc,\PP)},$ and similarly the
mapping telescope $\widetilde{\grO^2X}$ of the stabilized mapping space; some of the 
restriction maps are only
 homology fibrations, and the final result is that one has
isomorphisms
$$j_*:H_*(\widetilde{H^0(}\bbc,\PP))\ra{1.5} H_*(\widetilde{\grO^2X}).
\leqno{\stpp.13}$$ 
Using then the fact that the components of $\grO^2X$ are homotopy equivalent to
each other we have

\noindent{\sc Theorem} \stpp.14: \tensl Let $X$ be a smooth compact principal
almost solv-variety, and let $\bfk$ be a good multi-degree. Then there exists
a constant  $c_0 >0$ such that for all $t< c_0l(\bfk)| - 1$ there
are homology isomorphisms:
$$j_*: H_t(H^0_{\bf k}(\bbc,\PP))\ra{1.5} H_t(\grO^2_{\bfk}X).$$ \tenrm

We are now ready to extend the theorem to the more general case of theorem A;
this extension, in essence, adds to the action of our linear algebraic group
that of a compact complex Lie group (a torus). In terms of maps, this makes
essentially no difference, as the two-fold loop space of  a torus is
contractible;  in terms of holomorphic maps, it also makes no difference, 
since the holomorphic maps of $\bbp^1$ into the torus are
constants:

\noindent{\sc Proof of Theorem A}: By hypothesis the K\"ahler manifold $X$ is almost homogeneous
with respect to the connected complex solvable Lie group $S.$ Then by a theorem
of Remmert and Van de Ven [Akh] the Albanese map $\gra:X\ra{1.3} \hbox{Alb}(X)$
is a locally trivial holomorphic fibre bundle, and by a theorem of Oeljeklaus
[Akh], Oel] the fibre $F$ is a smooth connected, simply connected projective
algebraic variety. 
Moreover, the fibration is $S$-equivariant [HO,Huc], so
that $F$ is almost homogeneous with respect to a subgroup $S_0$ of $S,$ and has
an $S_0$-equivariant embedding of $F$ in some projective space. It follows that
$S_0$ is a solvable linear algebraic group and that $F$ is a smooth projective
algebraic variety with an open orbit biholomorphic to $S_0$ itself. In other
words $F$ is a smooth projective principal almost solv-variety. This allows us
to describe the space $\pi_2(X)$ of components of the mapping space
$\Omega_2(X)$ in the same way as before. Since the long exact sequence in
homotopy gives us an isomorphism  $\pi_2(F)\simeq \pi_2(X),$ we can then
think of $\pi_2(X)$ as being embedded in $H_2(F, N(F))= \bbz^r$, as before.  The
theorem follows by noticing that on the one hand, every map from $\bbp^1$ to a
torus is null-homotopic, and on the other, that every such holomorphic map is
constant.  Everything then lives in the fiber, and so it suffices to prove the
theorem for maps into $F.$ But this is the content of Theorem \stpp.14.
\hfill\za

We now examine how one can evaluate  the constant $c_0$. As we saw, the 
stability range of each stratum's contribution to the spectral sequence was 
governed by two numbers. 

a) The stability range for the homology in each stratum. This is  
$1/2$ of the minimum over the different components $X_\alpha$ 
of the number of simple $\alpha$-poles of elements of  the stratum. 

b) The codimension of the stratum within the space of holomorphic maps.
We saw that this was governed in turn by the codimension $\ell$ in the space
of $k$-jets of maps from $\bbc$ to $X$ of the subvariety of jets intersecting
$X_\infty$ with multiplicity $k$. 

Let us suppose that $X_\infty$ has smooth components with normal crossings. 
Then it is easy to see from expressions in local coordinates 
that the complex codimension $\ell$ is just $k$. From this, one then has that
the strata with $\rho$ simple poles in  $H^0_{{\bf k}}(\bbc,\PP)$ has complex
codimension $|{\bf k}|-\rho$, hence real codimension $2(|{\bf k}|-s)$.
Combining this with the stability range for the homology of the stratum gives
$c_0=1/2$.

To improve on this one must, in essence, compute some differentials in the
spectral sequence. It turns out that just computing one differential improves
things considerably.  Indeed, the cycles of dimension between $\rho/2$ and
$\rho$ which can cause the stability theorem to fail all contain factors (under
loop sum and bracket operations) of the form $[e,e']$, the one-dimensional
cycle  obtained by moving an $\alpha$-pole labelled by the base point $e$ in
$\LPP_{\alpha}$ in a circle around a fixed $\beta$-pole labelled by the base
point $e'$  in $\LPP_{\beta}$. One simply has to show that this cycle is killed
when one adds in the higher co-dimensional strata. To do this, one  adds in the
stratum with a double pole instead of the two simple poles. Geometrically, the
cycle, which is a circle, is filled in to a disk in which the two poles are
allowed to coalesce; this is quite easy to do in local coordinates, using the
normal crossing condition, and moving an embedded disk. The fact that these
cycles then do not contribute to the homology of the total space improves the
stability range to $c_0 = 1$, again when $X_\infty$ has smooth components with
normal crossings. A similar calculation, showing that a similar cycle is
trivial  in the case of $SU(2)$ instanton moduli, is given in more detail in
[Hu3].  It also increases the stability range there to the optimal $c_0=1$.

\bigskip

\noindent{\S5. \bf The Fundamental group and covering spaces}

\medskip

\noindent{\it (i) The fundamental group}

A crucial ingredient in the proof of the homotopy stability of our mapping
spaces is that $\holk$ has Abelian fundamental group when the components of
$\bfk$ are large enough. We have:

\noindent{\sc Theorem \ahpi.1:} \tensl If the good multi-degree $\bfk
=(k_1,\cdots,k_\kappa)$ satisfies $k_i\ge 2$ for all $i$, then the fundamental
group $\pi_1(\holk)$ is Abelian. \tenrm

\noindent{\sc Remark:} In the stable range, one then has that the fundamental group
of $\holk$ is simply the fundamental group of the mapping space, that is $\pi_3(X)$.

\noindent {\sc Proof:}  We begin by remarking that for dimensional
reasons (Theorem \stpp.3)  any loop in $H^0_\bfk(\bbc, \PP) = \holk$ can be deformed to the
generic  locus ${\cal S}$  consisting of configurations of points with simple
poles; this means that $\pi_1({\cal S})$ surjects onto $\pi_1(H^0_\bfk(\bbc, \PP))$.
  The locus ${\cal S}$ can be described as follows. One has $\kappa$ different
colors, corresponding to the different types of poles, that is the different
divisors at infinity, and associated to each color, a label space $\LPP_i$ of
principal parts of multiplicity 1 in the $i$-th color, 0 in the others.
The space $\LPP_i$ is connected, as it is a quotient of the space of  germs of 
maps into $X$ meeting the $i$-th component of infinity transversally, in its
smooth locus, and this latter space is connected. The generic stratum ${\cal S}$
then  consists of configurations of $|\bfk|=\sum_ik_i$ distinct
unordered points in the plane, with, for each $i$, $k_i$ points labelled by
elements of $\LPP_i$.  In other words, consider the product $\Pi_{i=1}^\kappa \bbd\bbp^{k_i}$
of the deleted symmetric products of $k_i$ points in $\bbc$ (i.e. unordered
$k_i$-tuples of distinct points in $\bbc$). Let ${\cal D} $ be the subspace of
$\Pi_{i=1}^\kappa \bbd\bbp^{k_i}$ consisting of those elements
$(\alpha_1,\cdots,\alpha_\kappa)$ which considered as $|\bfk|$ points of $\bbc$
have no points in common. One has a fibration $$ \Pi_i(\LPP_i)^{k_i}\rightarrow
{\cal S}\ra{1.5}{\cal D}.\leqno \ahpi.2$$ Let $G$ be the subgroup of the
braid group $B_k$  on $|\bfk|= \sum_ik_i$ letters consisting of those elements
which map to $S_{k_1}\times\cdots\times S_{k_\kappa}$ under the standard
homomorphism $B_{|\bfk|}\rightarrow S_{|\bfk|}$ of the braid group onto the symmetric group
$S_{|\bfk|}$ on ${|\bfk|}$ letters. ${\cal D}$ is a $K(G,1)$, and the $\LPP_i$ are connected,
so one has the exact sequence of groups
$$ 0\rightarrow\pi_1(\Pi_i(\LPP_i)^{k_i})\ra{1.5}\pi_1({\cal S})\rightarrow
\pi_1({\cal D})\ra{1.5}0.\leqno{\ahpi.3}$$ 
Using  the translation action on principal parts, one can lift paths in the base
${\cal D}$ to  ${\cal S}$. The monodromy of this lift along a loop $\ell$ in ${\cal D}$
is simply given by the action
of the image of $\ell$ in $S_{k_1}\times\cdots\times S_{k_\kappa}$ on
$\Pi_i(\LPP_i)^{k_i}$.  The fibration thus admits  sections: one simply
chooses the same label in $\LPP_i$ for each point of color $i$. This allows us
to split the sequence (\ahpi.3) and so decompose any element of $\pi_1({\cal
S})$ into a product of a loop in the fibers and a loop in $\pi_1 ({\cal D})$. 

The fundamental group of the full $\holk$ is a quotient of that of 
${\cal S}$ and one now would like to see what additional relations 
on $\pi_1({\cal S})$ are given by the inclusion of $\cals$ into $\holk$.
These relations are given by glueing in the complex codimension one strata.
The key relation turns out to be that particles of the same color are allowed
to coalesce: as in [BHMM1], section 9, one can bring two  principal parts of
multiplicity one and of the same color together into a double principal part.
 This kills the braiding between particles of the same color. To realise this coalescence 
for a given color $\rho$ in terms of a family of germs, one can choose
a coordinate $x$ on  $X$ such that $X_\rho$ is cut out locally by $x=0$. 
Now consider a family of curves $f_a: \bbc\rightarrow X$
whose $x$-component is given by $x(t,a) = t^2-a$; for $a\neq 0$, one has two
simple poles at $x=\pm\sqrt{a}$, and for $a=0$ one has a double pole.
Glueing in these extra strata includes ${\cal S}$ into  a bigger space ${\cal S}'$ 
for which double poles are allowed. Similarly,  ${\cal D}$ is glued into 
  a bigger space ${\cal D}'$ 
for which double points are allowed. Let $i$ denote this inclusion map. 

 The fundamental
group of ${\cal D}'$ is Abelian, if all the $k_i$ are at least two.
Indeed, let $\alpha$ be a standard braid group generator in $\pi_1({\cal D})$,
interchanging two poles of the same color, while leaving the others fixed; 
the preceeding paragraph tells us that $i_*(\alpha)=e$, the identity element; more generally,
$i_*(\alpha\beta\alpha^{-1}) = i_*( \beta)$. This means that one can use {\it any} of the
 points of a given color to represent in ${\cal D}$ an element of $\pi_1({\cal D}')$.
The natural thing is then to use only one of our points:
thus, if ${\cal D}_0$ is the space of $\kappa$-tuples of distinct points in the
 unit disk in $\bbc$,
one for each color, and if $j:{\cal D}_0\rightarrow {\cal D}$ is the stabilisation map of adding in 
fixed extra points outside the disk, 
then the composition $i\circ j$ induces a surjection on $\pi_1$.

Now let us suppose that there are two points $p_\ell$ and $q_\ell$ of each colour $\ell$  .
One can then represent an element $\gamma$ 
of $\pi_1({\cal D}')$ as a braid in the $p_\ell$ and another element $\delta$
as a disjoint (i.e., lying in another disk) braid in the $q_\ell$. In other words,
we can represent $\gamma,\delta$ as elements $(\hat\gamma, 0), (0,\hat\delta)$ in 
$\pi_1({\cal D}_0\times{\cal D}_0)$. 
Such elements commute, and so $\pi_1({\cal D}')$ is commutative.

The same holds true when one considers the space of labelled configurations: 
the key point is that the loop $\alpha'$ which interchanges two points 
with the same labels is contractible. This again allows us to ``localise''
loops onto configurations involving only one pole of each colour; if we have at our 
disposal at least two points of each colour, we can localise any two loops into loops 
of configurations living in disjoint disks, which must commute. In short, $\pi_1({\cal S}')$
is Abelian.

Thus, the relations put in by allowing
two points of the same color to coincide make the group Abelian, and as
$\pi_1(\holk)$ is obtained by putting in possibly more relations (adding more strata
and so going to an
even smaller quotient), it too must be Abelian. \hfill\za

\noindent{\sc Example \ahpi.4:}  The classical case of holomorphic maps of
the Riemann sphere to itself is a good example here.  The solvable group
that we select is $\bbc^*$ acting on $\bbp^1$, fixing zero and infinity.
In this case the open orbit is just a copy of $\bbc^*$, and ``$\infty$"
is $0 \sqcup \infty$.  Thus the data
describing a holomorphic map are
its zeros and poles.,
Any zeros can come together as can any poles, but roots and poles
cannot come together.  Thus, each stratum will consist of
pairs $\langle x_1, \dots, x_n\rangle$, $\langle y_1, \dots, y_m\rangle$
(the $x_i$'s are zeros and the $y_j$  poles). Attached to these points are
integers which are the multiplicities of the zeros and poles,
with the sole requirement that the sums of the multiplicities of  zeros and
the sums of the multiplicities of poles be equal to a fixed degree $n$.  
The fundamental group is  $\bbz$, with generator induced by looping 
a zero around a fixed pole [Se].
Note that the label spaces are discrete finite  sets and do not introduce
any generators into the fundamental group.

\noindent{\it (ii) Coverings of the strata}

We have just seen that the fundamental groups of our spaces of holomorphic maps 
are stably the first homology groups, and we know that these groups stabilize.
The stabilization maps extend to the universal cover; we have already seen that they 
preserved our strata.
Our next step is to see how the universal cover of $\holk$ behaves when one
restricts to the strata ${\cal V}_{\cal M}$. These strata correspond to
multiplicity patterns for the poles; there will be a certain number of {\it
simple} poles of each color, say $n_i$ for the $i$-th color, $i= 1,\dots,k$
and then multiple poles, of which there are, say $n_i, i = k+1,..., \ell$ of
various types $i$, taking values in principal part spaces which we also label
as $\LPP_i$. We restrict our attention to those strata with at least one simple
pole of each color, and note that the codimension of those strata with no
simple poles in a given color increases linearly with $|\bfk|$. Set $n =
\sum_{i=1}^\ell n_i.$ We can  describe the stratum ${\cal V}_{\cal M}$ as
$$C_n(\bbc)\times_{\Pi_iS_{n_i}}\Pi_i(\LPP_i)^{n_i}.$$ 
Here $C_j(X)$ denotes the configuration space of ordered $j$-tuples of distinct
points in $X$ and the symmetric group $S_j$, and its subgroup ${\Pi_iS_{n_i}}$
act by changing the ordering.

The restriction $\tilde{\cal V}_{\cal M}\rightarrow{\cal V}_{\cal M}$ of the
covering of $\holk$ is Abelian and is determined by a subgroup of  the first
homology group of ${\cal V}_{\cal M}= C_n(\bbc)\times_{\Pi_i{\cal
S}_{n_i}}\Pi_i(\LPP_i)^{n_i} $.  There is a fibration
$$\Pi_i(\LPP_i)^{n_i}\ra {1.5}C_n(\bbc)\times_{\Pi_i{\cal
S}_{n_i}}\Pi_i(\LPP_i)^{n_i} \ra{1.5}C_n(\bbc)/{\Pi_iS_{n_i}} \leqno \ahpi.6$$
which, as above, splits. As the base is an Eilenberg-Maclane space,
we have a semi-direct product 
$$\pi_1({\cal V}_{\cal M}) = [\Pi_i(\pi_1(\LPP_i))^{n_i}] \bowtie \pi_1
(C_n(\bbc)/{\Pi_iS_{n_i}}).$$ 
For the Abelianization of $\pi_1 (C_n(\bbc)/{\Pi_iS_{n_i}})$,
one sees that there are various winding numbers one can define: we are dealing
with $n_i$ particles for each color $i$, which can move around each other.  We
therefore have the total winding numbers for particles of color $i$ moving
around particles of color $i'$, where $i$ can equal $i'$, giving a factor of
$\bbz^{\ell(\ell+1)/2 }$ in the Abelianization. On the other hand, the
fundamental group is generated by elementary loops which send a particle of
color $i$ around a particle of color $i'$; the braiding relations tell us in
essence that when we Abelianize, we cannot tell which particle of color $i$ we
are sending around, so that the Abelian quotient is at most
$\bbz^{\ell(\ell+1)/2 }$.  In turn , then we have that $H_1({\cal V}_{\cal M})$
will be a quotient of $[\Pi_i (H_1(\LPP_i))^{n_i}] \times \bbz^{{\ell(\ell+1)/2
}}$, where we must still quotient out the action of the fundamental group of
the base (by conjugation) on that of the fibers; as this action is by
permutation of particles of the same color, one obtains for the first homology
group: $$H_1({\cal V}_{\cal M})= [\Pi_i (H_1(\LPP_i))] \times
\bbz^{\ell(\ell+1)/2 }$$ We note in particular that restricting any Abelian
covering of ${\cal V}_{\cal M}$ to the fibers of the fibration \ahpi.6 is
induced by addition maps $H_1(\LPP_i)^{n_i}\ra{1.5} H_1(\LPP_i)$.  Also, the
proof of theorem \ahpi.1 shows that our coverings quotient out the winding of
particles labelled by simple poles of the same color.

We can further decompose the space ${\cal V}_{\cal M}$ as an iterate
fibration, by successively projecting out the points labelled by $\LPP_\ell$,
then those labelled by  $\LPP_{\ell-1}$ and so on. In this way we obtain a
sequence of fibrations:
$$\eqalign{ 
{\cal V}_{\ell-1} ~\ra{1.5}~{\cal V}_{\cal M}~&\ra{1.5}
~C_{n_\ell}(\bbc)\times_{S_{n_\ell}}\LPP_\ell^{n_\ell} \cr
{\cal V}_{\ell-2}~\ra{1.5}~{\cal V}_{\ell-1}~&\ra{1.5}~C_{n_{\ell-1}}
(\bbc - \{y_{\ell,1},\cdots,y_{\ell,n_\ell}\} )
\times_{S_{n_{\ell-1}}}\LPP_{\ell-1}^{n_{\ell-1}}\cr
 {\cal V}_{\ell-3}~\ra{1.5} ~{\cal V}_{\ell-2}~&\ra{1.5}~C_{n_{\ell-2}}
(\bbc-\{ y_{\ell,1},\cdots,y_{\ell,n_\ell},
y_{{\ell-1},1},\cdots,y_{{\ell-1},n_{\ell-1}}\})\times_{S_{n_{\ell-2}}}
\LPP_{\ell-2}^{n_{\ell-2}}\cr
~\vdots~\phantom{\ra{1.5}} ~ \vdots
~&~\phantom{\ra{1.5}}\vdots\cr}\leqno\ahpi.7$$
$$
C_{n_1}
(\bbc- \{ y_{\ell,1},\cdots,y_{2,n_2}\})\times_{S_{n_1}} \LPP_1^{n_1}
~\ra{1.5}~{\cal V}_2~\ra{1.5}
C_{n_2} (\bbc-\{ y_{\ell,1},\cdots,y_{3,n_3}\})\times_{S_{n_2}}\LPP_2^{n_2}.
$$
All of these fibrations have sections. There is then a corresponding diagram
for the fundamental groups, which at each step is exact and splits as a
semi-direct product, and which at the i-th fibration projects out the  loops
involving the i-th particle space and the paths of particles of type $i$ around
fixed punctures corresponding to particles of type $i+1,\cdots,\ell$.
Restricting the maps of fundamental groups to the subgroup defining the covers,
we obtain a corresponding diagram of iterate fibrations for the covers
$\tilde{\cal V}_{\cal M}.$ Let $\tilde C_{n_i|n_{i+1},\cdots,n_\ell}(\LPP_i)$
denote the cover of $C_{n_i} (\bbc- \{
y_{\ell,1},\cdots,y_{i+1,n_{i+1}}\})\times_{S_{n_i}} \LPP_i^{n_i}.$ Then we
have the sequence of fibrations
$$\eqalign{ \tilde{\cal V}_{\ell-1} ~\ra{1.5}~\tilde{\cal V}_{\cal M}~&\ra{1.5}
\tilde C_{n_\ell|\cdot}(\LPP_\ell)\cr \tilde{\cal
V}_{\ell-2}~\ra{1.5}~\tilde{\cal V}_{\ell-1}~&\ra{1.5}~ \tilde
C_{n_{\ell-1}|n_\ell}(\LPP_{\ell-1})\cr
 \tilde{\cal V}_{\ell-3}~\ra{1.5} ~\tilde{\cal V}_{\ell-2}~&\ra{1.5}
\tilde C_{n_{\ell-2}|n_{\ell-1}, n_\ell}(\LPP_{\ell-2})\cr
~\vdots~\phantom{\ra{1.5}} ~ \vdots ~&~
\phantom{\ra{1.5}}\vdots\cr}\leqno\ahpi.8$$
$$\tilde C_{n_1|n_2,\cdots,n_\ell}(\LPP_1) ~\ra{1.5}~{\cal V}_2~\ra{1.5}
\tilde C_{n_2|n_3,\cdots,n_\ell}(\LPP_2).$$

\noindent{\sc Proposition} \ahpi.9: \tensl All the coverings in the fibrations
$\ahpi.7$ are Abelian. \tenrm

\noindent {\sc Proof}: Let us write one of the fibrations in \ahpi.7
and their coverings in \ahpi.8 as
$$\matrix {\tilde{\cal V}'\ &{\buildrel\tilde
a\over\longrightarrow}&\tilde{\cal V}\ & {\buildrel \tilde
b\over\longrightarrow}&\tilde{\cal C}\ \cr \decdnar{}\sigma'& &
\decdnar{}\sigma& & \decdnar{}\rho\cr {\cal V}'\ &{\buildrel
a\over\longrightarrow}&{\cal V}\ &{\buildrel
 b\over\longrightarrow}&
{\cal C}\ }$$
The downward maps all induce injections on the level of fundamental groups,
and the bottom row splits $\pi_1({\cal V})$ as a semi-direct product of $\pi_1(
{\cal V}')$ and $\pi_1({\cal C})$. The covering $\tilde C\rightarrow C$
corresponds to setting  $\rho_*(\pi_1(\tilde{\cal C})) =
b_*\sigma_*(\pi_1(\tilde{\cal V}))$. One must check that the commutator of
$\pi_1({\cal C})$ lies in $\rho_*(\pi_1(\tilde{\cal C}))$, and similarily that
the commutator of $\pi_1({\cal V}')$ lies in $\grs'_*(\pi_1(\tilde{\cal V}'))$.
A small diagram chase gives this from the fact that both these commutators lie
in the commutator of $\pi_1({\cal V}),$ which in turn lies in
$\grs_*(\pi_1(\tilde{\cal V}))$.  We see then that all the covers are governed
by subgroups of the first homology groups. \hfill\za

We are interested in the effect on these spaces of stabilization by adding in
simple poles; we will examine these one color at a time. To a configuration in 
the space $H^0_{\bf k}(\bbc,\PP)$, with poles located at $z_1,\cdots,z_n$ in
$\bbc$, we add a pole located at $\hbox{max}(|z_i|) +1$, labelled by the base
point in $\LPP_1$.  As noted in the previous section, this map preserves our
stratification, and we are now interested in the effect of the map on the level
of the covers. In the stable range, the covers are governed by the first
homology of $H^0_{\bf k}(\bbc,\PP)$, which is constant.  The stabilization preserves
all the spaces in the iterate fibrations \ahpi.7, \ahpi.8, and in fact only
changes  the last fibers of \ahpi.7, \ahpi.8, adding in an extra particle.  For
these last fibers, we have a commutative diagram: 
$$\matrix{\tilde
C_{n_1|n_2,\cdots,n_\ell}(\LPP_1)&\rightarrow& \tilde
C_{n_1+1|n_2,\cdots,n_\ell}(\LPP_1)\cr \decdnar{}& & \decdnar{}\cr
C_{n_1}(\bbc-\{y_{\ell,1},..,y_{2,n_2}\})\times_{{\cal
S}_{n_1}}\LPP_1^{n_1}&\rightarrow&
C_{n_1+1}(\bbc-\{y_{\ell,1},..,y_{2,n_2}\})\times_{{\cal
S}_{n_1+1}}\LPP_1^{n_1+1} }$$ 
where $n_1$ denotes the number of simple poles of the given color.

To show that the homology of the universal cover $\widetilde{H^0_{\bf k}(\bbc,
\PP)}$ of $H^0_{\bf k}(\bbc, \PP)$ stabilizes through an appropriate dimension
$q$, we proceed as in \stpp.11, using the spectral sequence which assembles the
homology of the different strata $\tilde{\cal V}_{\cal M}$, which reduces the
problem to showing that the homology of the stratum $\tilde{\cal V}_{\cal M}$
stabilizes in dimensions $q-(\hbox{codim}(\tilde{\cal V}_{\cal M}))$.  In turn
for these strata, one can use the spectral sequences for the fibrations
\ahpi.8, which, as the bases are not changed by the stabilisation, reduces the
stability question to showing that the coefficient systems $H_i(\tilde
C_{n_1|n_2,\cdots,n_\ell}(\LPP_1))$ stabilize in an appropriate range.  Our
theorem B, with the appropriate range, then will follow from

\noindent {\sc Theorem} \ahpi.10: \tensl Let $s$ be the rank of $\pi_3(X)$.
 The   stabilization map 
$$\tilde C_{n_1|n_2,\cdots,n_\ell}(\LPP_1)\ra{1.5}\tilde
C_{n_1+1|n_2,\cdots,n_\ell} (\LPP_1)$$
induces isomorphisms in homology groups $H_i$ for $i<n_1-s-1$. \tenrm

\noindent The proof of this theorem occupies the next section of the paper.

Using \ahpi.10, then, one obtains as in section 4.
 that the homology of  $\widetilde{H^0_{\bf k}(\bbc, \PP)}$
stabilises in the range $t\le c_0\ell({\bf k}) -s-2$. One then has to make sure that the 
space is simple in this range of dimensions, that is that the fundamental group acts 
trivially on the homology of the covering. This can be seen as follows.
Let ${\bf k}_0 = (1,1,...,1)$. We saw above that the fundamental group of 
${H^0_{{\bf k}_0}(\bbc, \PP)}$ surjects under the stabilisation map to that of 
${H^0_{\bf k}(\bbc, \PP)}$. Let an element $\tau$ in the fundamental group of ${H^0_{{\bf k}_0}(\bbc, \PP)}$
act by the deck transformation $T$ on the universal cover, and let the corresponding
element $\tau'$ of the fundamental group of ${H^0_{\bf k}(\bbc, \PP)}$ act by $T'$; we have
a commuting diagram on the universal covers:
$$\matrix{ 
\widetilde{H^0_{{\bf k}_0}(\bbc, \PP)}\times \widetilde{H^0_{\bf k-{\bf
k}_0}(\bbc, \PP)} &{\buildrel{\tilde I}\over {\longrightarrow }} &
\widetilde{H^0_{\bf k}(\bbc, \PP)}\cr (T,1)\decdnar{}& & T'\decdnar{}\cr
\widetilde{H^0_{{\bf k}_0}(\bbc, \PP)}\times \widetilde{H^0_{\bf k-{\bf
k}_0}(\bbc, \PP)} &{\buildrel{\tilde I}\over {\longrightarrow }} &
\widetilde{H^0_{\bf k}(\bbc, \PP)}},\eqno\ahpi.11$$
where $\tilde I$ represents the loop sum map.
Now let us suppose that we have a homology cycle   $A$ in
$\widetilde{H^0_{\bf k}(\bbc, \PP)}$ in the stable range for 
${\bf k-{\bf k}_0}$; we can represent it on the left hand side of \ahpi.11 as
a class ${\rm pt}\times A $. For obvious geometrical reasons, $(T,1)$ acts
trivially on ${\rm pt}\times A $; but then, referrring to the above diagram,
$T'_*(A) = A$, and so our spaces are
 simple in the appropriate range.

Once one has this simplicity, one  can then apply, e.g. the result of Hilton
and Roitberg ([HR], cor 3.4) to show that the homology of the universal covers
of $\holk$ stabilises to that of the universal covers of the mapping spaces
$\hbox{Map}_{\bfk}(\bbp^1,X)$; this then proves theorem B.  

\bigskip 

\def\cover{6}
\def\cals{{\cal S}}

\noindent{\S6. \bf Abelian covers of Labelled Particle Spaces on Punctured
Planes}

\medskip

\noindent{\it (i) Abelian covers of the label spaces}

To prove \ahpi.10, we want to understand the behaviour of the homology of the covers
$\tilde C_{n_1|n_2,..,n_\ell}(\LPP_1)$ of $C_{n_1}
(\bbc- \{ y_{\ell,1},..,y_{2,n_2}\})\times_{S_{n_1}} \LPP_1^{n_1}$
 as we stabilize.
We first must understand what type of cover we are dealing with. 
>From (\ahpi.6) and the discussion following it, we have that the 
subgroup of the first homology group of $C_{n_1}
(\bbc- \{ y_{\ell,1},..,y_{2,n_2}\})\times_{S_{n_1}} \LPP_1^{n_1}$ 
which governs the cover is a subgroup of $H_1(\LPP_1;\bbz)\times 
\bbz^{\ell-1}$. The $\bbz^{\ell-1}$ factor 
corresponds to the total winding of particles of color 1 
around the punctures corresponding 
to the other colors. We note that, as in theorem 5.1,  
any two   simple poles with labels in $\LPP_1$ are 
allowed to coalesce in $H^0_{\bf k}(\bbc, X)$, and so for our covers
there is no factor $\bbz$ corresponding to 
the self-winding of  particles of the same color. 

We begin with the case of no punctures; one ends up with this case
by   restricting our covers to 
the  subspace  $C_{n_1}(\bbc)\times_{S_{n_1}} \LPP_1^{n_1}$
corresponding to a (non-holomorphic!) embedding  $\bbc\rightarrow
\bbc- \{ y_{\ell,1},..,y_{2,n_2}\}$ . Having done this, 
the connected components of the cover have as  
group of deck transformations a quotient $\pi$ of $H_1(\LPP_1)$; restricting 
to a fiber  $\LPP_1^{n_1}$
the covers correspond to the composition of the    diagonal map 
$H_1(\LPP_1)^{n_1}\rightarrow H_1(\LPP_1)$ with projection to $\pi$.

In what follows, let $X$ denote $\LPP_1$, $k$ denote $n_1$ and $\tilde X$ the covering of $X$ 
with $\pi$ as deck transformations. Suppose that $\tilde X$ is given by a 
classifying map $\rho$:
$$\matrix{ \tilde X &\fract{\tilde \rho}{\ra{1.5}}& E_{\pi}\cr
\downarrow & &\downarrow\cr
X &\fract{\rho}{\ra{1.5}}& B_{\pi}\cr}\leqno \cover.1$$
where $\pi$ is an abelian group.  Then, the multiplication homomorphism
$\pi\times \pi \ra{1} \pi$, $(a,b) \mapsto ab$ induces a map
$$\grm_k\colon
\underbrace{B_\pi \times B_\pi \times \cdots \times B_\pi}_{k-times} \ra{1.5}
B_\pi$$
and consequently the composition $\grm_k\circ \rho^k\colon X^k \ra{1} B_\pi$.  
Moreover, $\grm_k\circ \rho^k$ is invariant under
the permutation action of the symmetric group, ${\cal S}_k$, on $X^k$.  
Consequently the induced $\pi$-covering map
$$\grp_k\colon {\tilde X}_k \ra{1.5} X^k\leqno{\cover.2}$$
is ${\cal S}_k$-equivariant, and the connected components of 
$\tilde C_{n_1|n_2,..,n_\ell}(\LPP_1)$ over  $C_k(\bbc)\times_{S_k} X^k$ are
given by $C_k(\bbc)\times_{S_k} \tilde X_k$. We note that ${\tilde X}_k =
(({\tilde X})^k)/\pi^{k-1}$ but $\grp^{k-1}\subset \grp^k$ is not the usual
inclusion and, in the end, this is what will cause us to have to take extra
care in our calculations.

\noindent{\it (ii) A spectral sequence for covers $C_k(\bbc)\times_{S_k}
\tilde X_k$} 
\medskip

Determining the cohomology of the resulting total space is naturally
somewhat involved.  In order to provide sufficient information to show
that we get stabilization we will have to introduce  somewhat novel
spectral sequences converging to the cohomology of 
 $\tilde X_k$ and that of the associated spaces 
$C_k(Y)\times _{{\cal S}_k}({\tilde X}_k)$
 which take as input the cohomology of $X$, and some data associated to the 
classifying map of $\tilde X$.

We shall work on the level of chain complexes. 
For $M$ a space, let $C_{\#}(M)$ be the  chain complex of $M$ with 
coefficients in a field $\bbf$. Set $E=E_\pi$. We note that $C_{\#}(\tilde X),
C_{\#}(E)$ can be chosen  so that the group ring $\bbf(\pi)$ acts freely on
them. 

\noindent{\sc Proposition \cover.3}:  \tensl Up to (free) equivariant chain
homotopy equivalance we can replace the (equivariant)
chain map $\rho_{\#}\colon C_{\#}({\tilde X})
\ra{1} C_{\#}(E)$ associated to the classifying map $\rho$ by
a surjection to $C_{\#}(E)$ of a free complex equivariantly homotopy equivalent
to $C_{\#}({\tilde X})$.\tenrm

\noindent(As this is standard we defer the proof to the end of this section.)

Then, if necessary replacing $C_{\#}({\tilde X})$ by the construction of the
proposition, we assume $C_{\#}({\tilde X})$ is given together with a surjection
$\rho_{\#}$ to $C_{\#}(E)$.  Consequently we have the short exact sequence of
($\bbf(\grp)$-free) chain complexes
$$ 0 \ra{1.5} L_{\#} \ra{1.5} C_{\#}({\tilde X})
\fract{\tilde\rho_{\#}}{\ra{1.5}} C_{\#}(E) \ra{1.5} 0.
\leqno{\cover.4}$$

\noindent{\sc Remark}:   Tensoring with $\bbf$ over
the group-ring $\bbf(\grp)$ and passing to cohomology gives $H^*(B_\pi,\bbf)$
 from the complex $C_{\#}(E)$, and $H^*(X,\bbf)$ from the complex 
 $C_{\#}({\tilde X})$.  The exact sequence above yields
 a long exact sequence, and the five-lemma then tells us that up to a shift,
$$H^*(B_{\grp}, X;\bbf) = 
H^*(L\otimes_{\bbf(\pi)}\bbf).\leqno{\cover.5}$$ 
We then have
$$\matrix{\cdots   \ra{1}&
 H^i(C_{\#}(\tilde X)\otimes_{\bbf(\pi)}\bbf)& \ra{1}&
H^{i}(L\otimes_{\bbf(\pi)}\bbf) &\ra{1}& H^{i+1}(C_{\#}(E)\otimes_{\bbf(\pi)}\bbf) &
 \cdots\cr
&\parallel&&\parallel&&\parallel\cr
\cdots   \fract{\rho^*}{\ra{1}}&
 H^i(X;\bbf)& \fract{j}{\ra{1}}&
H^{i+1}(B_{\grp}, X; \bbf) &\fract{\grd}{\ra{1}}& H^{i+1}(B_{\grp};\bbf) &
 \cdots.}\leqno{\cover.6}$$

>From  the $\bbf(\pi)$ chain subcomplex $L=L_{\#}$ of  $C_{\#}(\tilde X)$,
 we can build  a spectral sequence converging to 
$H^*(\tilde X_k, \bbf) = 
H^* (C_{\#}(\tilde X)^{\otimes k}\otimes _{\bbf(\pi^{k-1})}\bbf)$.
Let $A \subset H^*(X; \bbf)$ be the image
of $\rho^*\colon H^*(B_{\grp};\bbf) \ra{1} H^*(X; \bbf)$.  (For example,
if ${\tilde X} \ra{1} X$ is a $\bbz$-cover then it is classified
by a map $\rho\colon X \ra{1} B_{\bbz} \simeq S^1$, and $A = \langle
\rho^*(\gri)\rangle$ with $\gri$ the generator of $H^1(S^1)$.) The exact
sequence above relates  $A$
 to the relative groups $H^*(B_{\grp}, X;
\bbf)$ and will give rise shortly to the appearance of $A$ in the 
$E_2$-term of the spectral sequence.

\noindent{\sc Proposition \cover.7}:  \tensl
There is a spectral sequence converging to 
$H^*({\tilde X}_k; \bbf)$ with $E_1$-term
$$\eqalign{&\left[\bigoplus_{r=0}^{k-1}\left\{ H^*(B_{\grp}, X;\bbf)^{\otimes r}
\otimes H^*(B_{\grp};\bbf)^{\otimes k-r-1}\right\}
\otimes_{\bbf({\cal S}_r\times{\cal S}_{k-r})} \bbf({\cal S}_k))\right]\cr
&\hbox to 2.6in{\hfill}
\oplus \hbox to .1in{\hfill}H^*(B_{\grp}, X;\bbf)^{\otimes k}.\cr}
$$
and $E_2$-term
$$\left[\bigoplus_{r=1}^{k-1}\left\{ 
\left[{\tilde H}^*(X;\bbf)/A\right]^{\otimes r}\otimes A^{\otimes k-r-1}
\right\}\otimes_{
\bbf({\cal S}_r\times{\cal S}_{k-r})}\bbf({\cal S}_k)\right]
\hbox to .1in{\hfill}
\oplus \hbox to .1in{\hfill}L_k$$
where both $H^*(B_{\grp}, X;\bbf)^{\otimes k}$ and $L_k$ are at least $k$-connected.\tenrm

\noindent{\sc Remark}:  The action of ${\cal S}_i$ on
$H^*(B_{\grp}, X;\bbf)^{\otimes i}$ is just permutation of coordinates
while that of ${\cal S}_{k-i}$ on $H^*(B_{\grp};\bbf)^{\otimes k-i-1}$ 
is obtained from an ${\cal S}_{k-i}$-equivariant embedding of 
$B_\pi^{k-i-1}$ into $B_\pi^{k-i}$ derived from  the embedding $\grp^{k-i-1}
\subset \grp^{k-i}$ as the subgroup of all elements $(g_1, \dots,
g_{k-i})$ with product $g_1g_2\cdots g_{k-i} = 1$.  The fact that we
are not dealing with a standard permutation action here will be the major
difficulty with proving stability.

\noindent{\sc Proof}:
The inclusion $L_{\#} \subset C_{\#}({\tilde X})$
gives rise to an ${\cal S}_k$-equivariant
filtration of the $k$-fold tensor product, $( C_{\#}({\tilde X}))^k$:
$$L^k \subset (L^{k-1}\otimes C + L^{k-2}\otimes C\otimes L + \cdots + C\otimes L^{k-1}) \subset
(L^{k-2}\otimes C^2 + \cdots + C^2\otimes L^{k-2}) \subset \cdots \subset C^k$$
where $L= L_{\#}, C = C_{\#}({\tilde X})$, and the exponents denote tensor product. 

The relative quotients are of the form
$$L^k, (L^{k-1}\otimes C(E) \oplus \cdots \oplus C(E)\otimes L^{k-1}),~(L^{k-2}\otimes
C(E)^2 \oplus \cdots \oplus C(E)^2\otimes L^{k-2}),\cdots$$
and taking $\otimes_{\bbf(\pi^{k-1})}\bbf$,  we have an  associated spectral sequence for ${\tilde X}_k =
(\tilde X^k)/\pi^{k-1}$.

Now we show that  the $E_1$-term of this spectral
sequence is of the form desired. As the tensoring with $\bbf({\cal S}_k)$ in the statement of the proposition 
 simply accounts for the different possible 
orderings of the factors, it will suffice to show that
 
$$\eqalign{H^*((L^{\otimes r}\otimes C_{\#}(E)^{\otimes s})
\otimes_{\bbf(\grp^{r+s-1})}\bbf)&=
 H^*(L\otimes_{\bbf(\grp)}\bbf)^{\otimes r }\otimes
H^*(C_{\#}(E)\otimes_{\bbf(\grp)}\bbf)^{\otimes s-1}\cr
 & =H^*(B_\pi,X;\bbf)^{\otimes r }\otimes
H^*(B_\grp;\bbf)^{\otimes s-1}},\leqno \cover.8 $$
where $r+s=k$, and the action of $(g_1,..,g_{r+s-1}) \in \grp^{r+s-1}$ is via
the rule
$$g_i(x_1\otimes \cdots \otimes x_{r+s})
= (x_1\otimes \cdots \otimes g_i(x_i)\otimes \cdots \otimes g_i^{-1}(x_{r+s})).
\leqno{\cover.9}$$
By this we mean that $g_i$ only acts non-trivially on the $i^{th}$
and $k^{th}$ coordinates and there acts as specified. Taking the natural
filtration by dimension on the first $k-1$ factors gives (yet another)   spectral sequence
whose $E_2$-term is the cohomology of the complex
$\hbox{Hom}_{\bbf(\grp^{r+s-1})}(L^{\otimes r }\otimes C_{\#}(E)^{\otimes s-1},
H^*(C_{\#}(E);\bbf))$ , namely, recalling \cover.3,
$$H^*_{\grp^{r+s-1}}(L^{\otimes r }\otimes C_{\#}(E)^{\otimes s-1};
H^*(C_{\#}(E);\bbf)) =H^*((L\otimes_{\bbf(\pi)}\bbf)^{\otimes r }\otimes
(C_{\#}(E)\otimes _{\bbf(\grp)}\bbf)^{\otimes s-1})$$ concentrated along the
line $E_{*, 0}$ since the action of $\grp$ on
$$H^*(C_{\#}(E);\bbf) = \cases{\bbf& if $*=0$,\cr
0& otherwise\cr}$$
is trivial.  It follows that $E_2 = E_{\infty}$ for this sequence, giving the
equality \cover.8.

Now that we have evaluated the $E_1$-term of our spectral sequence, note that
there is a non-trivial $d_1$-differential associated to the differential in
the long exact sequence in cohomology of $(\cover.4)$ above.  The effect of
this differential, in view of the results above is to replace
$H^*((L\otimes_{\bbf(\pi)}\bbf)$ by $H^*(X)/A$ and $H^*(B_\grp;\bbf)$ by $A$
(see the remark following $(\cover.4)$).  Thus we have determined the
$E_2$-term of the spectral sequence, except for the part $L_k$ corresponding to
$H^*(L^{\otimes k}\otimes_ {\bbf(\grp^{k-1})}\bbf)$, which starts in dimension
at least $k$.

For later use we record the following consequence of \cover.8:

\noindent{\sc Lemma \cover.10}:  \tensl Recalling the identifications
\cover.6 and \cover.8, the natural map
$$H^*(B_{\grp}, X;\bbf)^{\otimes r }\otimes H^*(\grp;\bbf)^{\otimes s}
\ra{1} H^*((L^{\otimes r }\otimes C_{\#}(E)^{\otimes
s})\otimes_{\bbf(\grp^{r+s-1})}\bbf)$$
induced by the inclusion $\grp^{r+s-1}\hookrightarrow \grp^{r+s}$
of \cover.9 is a surjection and is
equivariant with respect to the action of ${\cal S}_r\times {\cal S}_s$
which permutes coordinates.\tenrm

Next, we have

\noindent{\sc Lemma \cover.11}:  \tensl Let $\grp = \prod_1^m\bbz/m_j \oplus
\bbz^t$ with generators $\gra_1, \dots, \gra_m$, $\grc_1,\dots, \grc_t$.  Let
$a_i \in H^1(B_\grp^s;\bbf)$, $i= 1, \dots ,m $ be $\sum_{j=1}^s \gra_{i,j}^*$,
$c_i$ in $H^1(B_\grp^s;\bbf)$ be the sum of elements dual to the integral
generators $c_i = \sum_{j=1}^s\grc_{i,j}^*$, $i = 1, \dots, t$, and $b_i =
\grb(a_i)$, the appropriate Bochstein of $a_i$.  Then, thinking of
$H^*(B_{\grp}, X;\bbf)^{\otimes r} \otimes H^*(B_\grp;\bbf)^{\otimes s}$ as
being embedded in the tensor algebra of $H^*(B_{\grp}, X;\bbf) \oplus
H^*(B_\grp;\bbf)$,  the kernel of the surjection above is the intersection of
$H^*(B_{\grp}, X;\bbf)^{\otimes r} \otimes H^*(B_\grp;\bbf)^{\otimes s}$ with
the ideal generated by the elements $a_i$, $b_i$, and $c_i$.  \tenrm

\noindent{\sc Proof}:  The point is to understand the inclusion of
$\grp^{s-1}\subset \grp^s$.  Note that by definition $\grp^{s-1}$ is the kernel
of the {\it sum} map: $\grm_s\colon \grp^s \ra{1} \grp$, so we have an exact
sequence $$\grp^{s-1} \fract{i}{\hookrightarrow} \grp^s
\fract{\grm_s}{\ra{1.5}} \grp$$ and the kernel of $i^*$ is evidently the ideal
generated by the image of $(\grm_s)^*,$ but the elements above are precisely
those images.\hfill\za

\noindent{\sc Proposition \cover.12}:  \tensl There is a spectral sequence
converging to $H^*(C_k(Y)\times _{{\cal S}_k}{\tilde X}_k)$ with the
$\grp^{k-1}$-action as above, that has $E_2$-term
$$E_2^i ~=~H^*_{{\cal S}_k}(C_k(Y); H^*(B_{\grp}, X;\bbf)^{\otimes i}\otimes
H^*(B_{\grp};\bbf)^{\otimes k-i-1}\otimes_{\bbf({\cal S}_i\times
{\cal S}_{k-i})}\bbf({\cal S}_k))$$
for $i< k$, while $E_2^k$ starts in dimension at least $k$
and $E_2^j = 0$ for $j > k$.\tenrm

\noindent{\sc Remark}: The actions of the symmetric groups on 
$H^*(B_{\grp}, X;\bbf)^{\otimes i}$ and 
$H^*(B_{\grp};\bbf)^{\otimes k-i-1}$  are the same as above, Also,
 since the ${\cal S}_k$-module above is
induced up from an ${\cal S}_i\times {\cal S}_{k-i}$-module, we
see that
$$\eqalign{H^*_{{\cal S}_k}(C_k(Y); H^*(B_{\grp},& X;\bbf)^{\otimes i}
\otimes H^*(B_{\grp};\bbf)^{\otimes k-i-1}\otimes_{\bbf({\cal S}_i
\times {\cal S}_{k-i})}\bbf({\cal S}_k))\cr
& \cong
H^*_{{\cal S}_i\otimes {\cal S}_{k-i}}(C_k(Y);
H^*(B_{\grp},X;\bbf)^{\otimes i}
\otimes H^*(B_{\grp};\bbf)^{\otimes k-i-1}).\cr}$$

\noindent{\sc Proof}:   
The filtration of $C_{\#}({\tilde X}^k)\otimes_{\bbf(\grp^{k-1})}\bbf$
constructed above is invariant under the action of ${\cal S}_k$.  Consequently,
it gives rise to a filtration of the chain complex for
$C_k(Y)\times_{{\cal S}_k}{\tilde X}^k/(\grp^{k-1})=
C_k(Y)\times_{{\cal S}_k}{\tilde X}_k$, and the $E_0$-level term for
this spectral sequence is \break
$Hom_{\bbf({\cal S}_k)}(C_{\#}(C_k(Y), E_0(fiber))$ with differential
dual to $\partial\otimes 1 + \gre (1\otimes d_0)$ where $d_0$ is the
first differential in the spectral sequence of the fiber.  
However, $C_k(Y)$ is a free $\bbf({\cal S}_k)$-module so the usual
spectral sequence for determining the cohomology of this complex
reduces to the sequence of its edge terms which gives
$E_1$ as a chain complex with (local) coefficients with values in the 
$E_1$ term of the fiber, and so $E_2 = H^*(C_k(Y); E_1(fiber))$ as
asserted.\hfill\za

\noindent{\sc Proof of Proposition \cover.3}:  Let $f\colon V_{\#} \ra{1}
W_{\#}$ be a chain map, then the {\tensl mapping cylinder} of $f$ is the
complex
$$MCyl(f)_{\#} = V_{\#} \oplus V_{\#-1} \oplus W_{\#}$$
with boundary defined by
$$\partial(v | v^{\prime} | w) = (\partial v + v^{\prime}|
-\partial v^{\prime} | -f(v^{\prime}) + \partial w).$$
It is chain equivalent to $W_{\#}$ via the map $\psi$:
$$\psi\colon (v|v^{\prime}|w) = w + f(v),$$
and converts the map $f$ into the inclusion $I\colon
V_{\#} \hookrightarrow MCyl(f)_{\#}$ with $I(v) = (v,0,0)$.
Similarly, the {\tensl mapping cone} of $f$ is
the chain complex $MC(f) = V_{\#-1}\oplus W_{\#}$
with boundary given as $\partial(v|w) = (-\partial(v)|\partial(w) + f_*(v))$.

In our case, in order to convert the chain map $\rho_{\#}$ into a surjection
we combine the two above constructions and replace $C_{\#}({\tilde X})$ by the
chain complex
$$C_{\#}({\tilde X})\oplus C_{\# + 1}(E_\pi)\oplus C_{\# -1}({\tilde X})\oplus
C_{\#}(E_\pi)\oplus C_{\#}({\tilde X})\leqno \cover.13$$
If we write an element as
$$\{c_n, e_{n+1} | c^{\prime}_{n-1}, e^{\prime}_n | c_n^{\prime\prime}\},$$
the boundary is defined by
$$\eqalign{\partial \{ (c_n,e_{n+1}&|c_{n-1}^{\prime}, e_n^{\prime}|c_n^{\prime\prime}\}
= \cr
&\{-\partial c_n + c_{n-1}^{\prime}, \rho_{\#}(c_n) + \partial e_{n+1} + e_n^{\prime}|
\partial c_{n-1}^{\prime}, - \rho_{\#}(c_{n-1}^{\prime}) - \partial e_n^{\prime}|
-c_{n-1}^{\prime} - \partial c_n^{\prime\prime}\}.\cr}$$
Now, if $(V_{\#}, \partial)$
is a chain complex, then the {\tensl suspension}
$(S(V)_{\#}, \partial)$ is defined as $(V_{\#},
-\partial)$.  The map $\gre\colon V_{\#} \ra{1} S(V)_{\#}$ defined by
$\gre(v) = (-1)^{|v|}v$ is easily seen to be a chain equivalence.
( Normally, the suspension $SV_{\#}$ is regraded, so that
$SV_n = V_{n-1}$, but this would introduce needless indices into our constructions,
hence we leave the grading degree as it was.  The important point for us is the
fact that the resulting complex is equivariantly chain equivalent to the
original one.) Here, the chain equivalence of (\cover.13) with $SC_{\#}({\tilde
X})$ is given by sending the $5$-tuple above to $c_n^{\prime\prime} + c_n$,
while the surjection to $SC_{\#}(E)$ is given by sending the $5$-tuple to
$e_n^{\prime} - \rho_{\#}{\grp}(c_n^{\prime\prime})$.  \hfill\za

\noindent{\it (iii) The cohomology of $C_k(\bbc)$ with coefficients
in the $E_1$-term of the fiber}
\medskip
We want the stabilization map on the last fiber of our iterate fibration,
$$\tilde C_{n_1|n_2,...,n_\ell}(\LPP_1)\rightarrow \tilde
C_{n_1+1|n_2,...,n_\ell}(\LPP_1),$$
to induce isomorphisms in homology through a suitable range. We first examine
the case $n_2=n_3=...n_\ell=0$, that is, when there are no punctures.
(Alternatively, one can look at the subspace of particles over $\bbc\subset
\bbc- \{ y_{\ell,1},..,y_{2,n_2}\}$, and take a connected component of the
cover.) Setting $n_1= k, \LPP_1=X$ as above, one  then has a map
$$C_k(\bbc)\times _{{\cal S}_k}({\tilde X}_k)\ra {1.5}
C_{k+1}(\bbc)\times _{{\cal S}_{k+1}}({\tilde X}_{k+1}).\leqno \cover.14$$

\noindent{\sc Proposition} \cover.15: \tensl The map \cover.14 induces
isomorphisms in cohomology in dimensions less than $k-s-1$, where $s$ is the
rank of $\pi.$ \tenrm

\noindent {\sc Proof}: In view of the above spectral sequences, 
we want to understand the stability properties of 
$$\eqalign{H^*_{{\cal S}_k}&(C_k(\bbc); H^*(B_{\grp}, X; \bbf)^{\otimes
i}\otimes H^*(B_{\grp};\bbf)^{\otimes k-i-1}\otimes_{\bbf({\cal S}_i\times
{\cal S}_{k-i})} \bbf({\cal S}_k))\ra{1.5}\cr &H^*_{{\cal
S}_{k+1}}(C_{k+1}(\bbc); H^*(B_{\grp}, X; \bbf)^{\otimes i}\otimes
H^*(B_{\grp};\bbf)^{\otimes k-i}\otimes_{\bbf({\cal S}_i\times {\cal
S}_{k-i+1})} \bbf({\cal S}_{k+1}))}\leqno{\cover.16}$$
We note that the stabilization adds a point labelled by the base point in
$X$, and so on the algebraic level is just tensoring by  a $0$-cycle.  Since
$H^0(B_{\grp}, X; \bbf)=0$, the stabilization factors through an increase in the
number of $H^*(B_{\grp};\bbf)$-factors. The space $C_k(\bbc)/{\cal S}_k$ is the
classifying space of the braid group $\Gamma_k$, which maps to the symmetric
group ${\cal S}_k$. We are therefore determining the stability properties of
the cohomology of the braid group with coefficients in the $\bbf({\cal
S}_k)$-modules above. Loop space theory tells us how to determine the
cohomology of $\grG_k$ with coefficients in modules of the form $A^{\otimes
k}$, where ${\cal S}_k$ acts to permute coordinates.   But the modules in
question here are not of this form.  So we will have to introduce further
arguments. We begin by giving the stability results for $H^*(C_k(\bbc);
A^{\otimes k})$ where $A$ is a tensor product $E(e_1, \dots, e_s)\otimes
\bbf[b_1, \dots, b_r]$ with $\dim(e_i) = 1$, $\dim(b_j) = 2$.

\noindent{\sc Lemma} \cover.17: ([BHMM2, lemma 6.8])  \tensl Let $A$ be as
above, then the cohomology map induced by the natural map of chain complexes
$${\cal C}_{\#}(C_k(\bbc)\otimes_{\bbf({\cal S}_k)}A^k)
\ra{1.5} {\cal C}_{\#}(C_{k+1}(\bbc)\otimes_{\bbf({\cal S}_k\times
1)}A^k\otimes 1) \ra{1.5} {\cal C}_{\#}(C_{k+1}(\bbc)\otimes_{\bbf({\cal
S}_{k+1})}A^{k+1})$$
is an isomorphism through dimensions less than $2k-s+1 $.\tenrm

We next determine stability dimensions for the map \cover.16
for the case where $\grp$ is {\it cyclic}:  in the cohomology of $\Gamma_k$, we
replace the tensor products $A^{\otimes k}$ by the modules in Lemma \cover.11.
Assume that $\grp = \bbz$ or $\bbz_{p^i}=\bbz/p^i\bbz$ with $p$ a prime, so
$H^*(\grp ;\bbf) = E(e_1)$ or $E(\gra_1)\otimes \bbf[\grb_2]$ and
$$H^*(\grp^k;\bbf) = \cases{E(\gra_1, \dots, \gra_k) & if $\grp  = \bbz$,\cr
E(\gra_1,\dots, \gra_k)\otimes \bbf[\grb_1,\dots, \grb_k]& if $\grp =
\bbz_{p^i}$.\cr}$$
We only do the case $\pi = \bbz,$ as the others are similar. We first determine
the structure of the ideal  described in Lemma \cover.11 and the corresponding
quotient.  We define $V(m) = E(\gra_1, \dots, \gra_m)/(\sum \gra_i)$ in the
notation of Lemma \cover.11 and we have

\noindent{\sc Lemma \cover.18}:  \tensl  There is a short exact sequence
$$0 \ra{1.5} V(m) \fract{\cup \sum \gra_i}{\ra{1.5}}
 E(\gra_1,\dots, \gra_m) \fract{\grp}{\ra{1.5}}
V(m) \ra{1.5} 0.$$
\tenrm

\noindent{\sc Proof}:   For convenience denote by $a(m)$ the sum $\Sigma_i\gra_i \in E(\gra_1,
\dots \gra_m)$. The kernel of the projection is composed of elements of the form
$\omega a(m)$, and there is a natural map from $V(m)$ to the kernel given
by $\omega\mapsto \omega a(m)$. It is obviously surjective, and to see
that it is injective, suppose that $\omega a(m) =0$. 
 Note that $a(m) = a(m-1) + \gra_m$, and we can write
$$E(\gra_1, \dots, \gra_m) = E(\gra_1, \dots, \gra_{m-1})
\oplus E(\gra_1, \dots, \gra_{m-1})\gra_m.$$
Set $\omega =\grT + \grn \gra_m$, then
$0= \omega a(m) = \grT a(m-1) + (- \grn a(m-1)+ \grT)\gra_m$, and so
$\grT = \grn a(m-1)$ giving $\omega = \grn a(m)$.\hfill\za

Now expanding out the exact sequence of Lemma \cover.18 with respect to the
grading gives
$$\eqalign{0 \ra{1.5} \bbf \ra{1.5}& E^1(\gra_1,\dots, \gra_m) \ra{1.5}
V^1(m)\ra{1.5} 0\cr 0 \ra{1.5} V^1(m) \ra{1.5}& E^2(\gra_1,\dots, \gra_m)
\ra{1.5} V^2(m) \ra{1.5} 0\cr 0 \ra{1.5} V^2(m) \ra{1.5}& E^3(\gra_1, \dots,
\gra_m) \ra{1.5} V^3(m) \ra{1.5} 0\cr \phantom{\vdots} \hbox to.5in{\hfill}
\vdots \hbox to.3in{\hfill}& \hbox to.5in{\hfill}\vdots \hbox to
.7in{\hfill}\vdots\cr}.\leqno \cover.19$$
>From Lemmas \cover.10 and \cover.11, we then have an exact sequence
$$\eqalign{ 0\ra{1.5} H^*(B_{\grp}, X;\bbf)^{\otimes r }\otimes
V^*(k-r)\ra{1.5}& H^*(B_{\grp}, X;\bbf)^{\otimes r }\otimes
H^*(\grp;\bbf)^{\otimes k-r}\cr &\ra{1.5} H^*((L^{\otimes r }\otimes
C_{\#}(E)^{\otimes k-r})\otimes_{\bbf(\grp^{k-1})}\bbf)\ra{1.5} 0.}$$
The stability result for homology with coefficients in the last module
will follow from that of the other two. The middle one has the stability
property in dimensions $<2k$, by Lemma \cover.18. For the first, we can build
up inductively using  \cover.19 and Lemma \cover.17.

We note that stabilization gives maps
$$\eqalign{H^*_{{\cal S}_k}(C_k(\bbc);
& H^*(B_{\grp}, X)^r\otimes V^*(k-r) \otimes_{\bbf({\cal S}_{r}\times {\cal
S}_{k-r})}\bbf({\cal S}_k))\cr &\ra{1.5} H_{\cals_{k+1}}^*(C_{k+1}(\bbc);
H^*(B_{\grp}, X)^r\otimes V^*(k-r+1)\otimes_{\bbf({\cal S}_r\times {\cal
S}_{k-r+1})}\bbf({\cal S}_{k+1})),\cr}$$
since $1 \not\in H^*(B_{\grp}, X)$. This fact also shifts the stability up by
$r$ dimensions. On the other hand each of the short exact sequences above gives
rise to a long exact sequence 
$$\eqalign{\cdots \ra{1.5} H^*_{{\cal S}_k}(&C_k(\bbc);
H^*(B_{\grp},X)^{\otimes r}\otimes V^i(k-r)\otimes_{\bbf({\cal
S}_{r}\times {\cal S}_{k-r})}\bbf({\cal S}_k))\cr  &\ra{1.5} H^*_{{\cal
S}_k}(C_k(\bbc);  H^*(B_{\grp}, X)^{\otimes r}\otimes E^{i+1}(\gra_1, \dots,
\gra_k)\otimes_{\bbf({\cal S}_{r}\times {\cal S}_{k-r})}\bbf({\cal S}_k))\cr
&\hbox to .4in{\hfill} \ra{1.5} H^*_{{\cal S}_k}(C_k(\bbc); H^*(B_{\grp},
X)^{\otimes r}\otimes V^i(k-r)\otimes_{\bbf({\cal S}_{r}\times {\cal
S}_{k-r})}\bbf({\cal S}_k))\ra{1.5} \cdots\cr}$$
which, in turn gives an inductive method for determining
$$\eqalign{H^*_{{\cal S}_k}(C_k(\bbc); 
 &H^*(B_{\grp},X;\bbf)^{\otimes r}\otimes V^i(k-r)\otimes_{\bbf({\cal S}_r
\times {\cal S}_{k-r})}\bbf({\cal S}_k))\cr
&= H^*(\grG_{k,r};H^*(B_{\grp},X;\bbf)^{\otimes r}\otimes V^i(k-r)),\cr}$$
where $\grG_{k,r}$ is the subgroup of $\grG_k$ given as the inverse image of
${\cal S}_{k-r} \times {\cal S}_r$ under the projection $\grG_k \ra{1} {\cal
S}_k$.  Of course, we don't need the entire calculation, but what is clear is
that the five lemma shows that we lose at most one dimension of stabilization
passing from $V^i(k-r)$ to $V^{i+1}(k-r)$.  But since, at the next stage the
dimension of $V^i(k-r)$ is augmented by one (since there it is $V^i(k-r)\cup
a(k)$), we actually do not lose this dimension. We then find that homology
stabilizes in the first term for dimensions less than $2(k-r) + r$  and
so the homology of the third term stabilizes in dimensions $2(k-r) +r -1\ge
k-1$.

There remains the question of stability for arbitrary $\pi$. 
The resulting stability dimensions will then decrease as we replace the cyclic
group by a general finitely generated abelian group by a simple function of the
number of generators of $\grp$ since the construction preserves tensor
product.  That is to say, if $\grp = \grp_1\times \grp_2$ then $H^*(\grp^k) =
H^*(\grp_1^k)\otimes H^*(\grp_2^k)$; similarly  the quotients $ H^*((L^{\otimes
r }\otimes C_{\#}(E)^{\otimes k-r}) \otimes_{\bbf(\grp^{k-1})}\bbf)$ for $\pi$
decompose into tensor products of the corresponding factors for $\pi_1$,
$\pi_2$. One then obtains stability for $\pi$-covers from that of the
$\pi_i$-covers, with the possible loss of one dimension due to ``edge
effects''. One then has the stability for the terms of the spectral sequence
\cover.16, and so, with the loss of one dimension, for the cohomology of the
spaces \cover.14. This concludes the proof of Proposition \cover.15. \hfill\za


\noindent {\it (iv) Putting in the punctures.}
\medskip

We exploit the  result for particles over $\bbc$ to obtain the corresponding
stability result for the strata $\tilde C_{n_1|n_2,...,n_\ell}(\LPP_1)$.
Keeping $n_i=k$, $X=\LPP_1$ as before, we renumber the punctures $y_{ij}$ as
$y_s, s=1, ..,m$, where $m= \sum_{i=2}^\ell n_i$,   and suppose that the points
$y_i$ are placed along the $y$ axis. Let $\bbr^+_i$ denote the half-lines
$\{(x,y_i)|x>0\},$ and denote by $D$ the complement of the $\bbr^+_i$. We
stratify $C_k(D-\{y_1,..,y_m\})\times_{{\cal S}_k}X^k$ according to the
number of particles on the  $\bbr^+_i$ with strata
$${\cal D}_{j, r_1,...,r_m} = (C_j(D)\times C_{r_1}\times..\times C_{r_m})
\times_{{\cal S}_j\times{\cal S}_{r_1}\times..\times {\cal S}_{r_m}}X^k.$$
As there is a natural ordering of points along
lines, one can rewrite this as $${\cal D}_{j, r_1,...,r_m} = (C_j(D)\times
\Delta^{k-j})\times _{{\cal S}_j} X^k$$ 
where $\Delta^{k-j}$ is an open $k-j$-cell, and ${\cal S}_j$ acts on
$C_j(D)$ and the first $j$ factors of $X^k$.
Passing to the covers, we have strata
$$\tilde {\cal D}_{j, r_1,...,r_m} = 
[(C_j(D)\times \Delta^{k-j})\times _{{\cal S}_j}\tilde X_k]\times G,
\leqno \cover.20$$
where $G$ is a discrete Abelian group, corresponding to pure winding around the
punctures. This factor appears because staying on the stratum precludes winding
around the punctures. We have the following easy lemma, giving the codimension
of the strata:

\noindent{\sc Lemma \cover.21}: \tensl The strata $\tilde{\cal D}_{j,
r_1,...,r_m}$ are smooth and  have real codimension $k-j$.\tenrm

We can compute the cohomology of $\tilde C_{n_1|n_2,...,n_\ell}(\LPP_1)$
from the spectral sequence associated to this stratification, in which 
the cohomology of the strata appears suspended by the codimension. 
The stabilization preserves the stratification and so the stability 
result for $\tilde C_{n_1|n_2,...,n_\ell}(\LPP_1)$ follows from the 
following:

\noindent{\sc Proposition   \cover.22}: \tensl The stabilization map
$$\tilde {\cal D}_{j, r_1,...,r_m}\rightarrow \tilde {\cal D}_{j+1,
r_1,...,r_m}$$
induces isomorphisms in cohomology in dimensions $<j-s-1$.\tenrm

\noindent{\sc Proof}: In analogy with \cover.12,
one has a spectral sequence for the strata  whose $E_2$ term is 
the cohomology of $C_j(Y)$ with coefficients in the ${\cal S}_j$-module
$$ H^*(B_{\grp}, X;\bbf)^{\otimes i}\otimes
H^*(B_{\grp};\bbf)^{\otimes k-i-1}\otimes_{\bbf({\cal S}_i\times
{\cal S}_{k-i})}\bbf({\cal S}_k)+ H_k,\leqno \cover.23$$
where $H_k$ is $k$-connected, and so doesn't enter into consideration.
To see what this  ${\cal S}_j$-module is, let us write 
$$\eqalign{ A=& H^*(B_{\grp}, X;\bbf),\cr B=& H^*(B_{\grp};\bbf).}$$
As a vector space, \cover.23 is 
$$\oplus_{(i, k-i){\rm shuffles}} A^{\otimes i}\otimes
B^{\otimes k-i-1}.$$
Think of each of these terms corresponding to the choice of $i$ ``slots''amongst
$k$ boxes; the action of ${\cal S}_j$ on these depends on the number $\ell$ of
these slots lying in the first $j$ boxes. As an ${\cal S}_j$-module, one can
identify it with
$$\oplus_{\ell=0}^{min(j,i)}\left(A^{\otimes \ell}\otimes
B^{\otimes j-\ell-1}\otimes_{\bbf({\cal S}_\ell\times {\cal
S}_{j-\ell})}\bbf({\cal S}_j)\right) \otimes \left((A^{\otimes i-l}\otimes
B^{\otimes k-i-j+\ell} \right)^{C(k-j,k-\ell)},$$ 
where $C( k-\ell, k-j)= (k-\ell)!/[(k-j)!(j-\ell)!]$.  The action of the ${\cal
S}_j$ on the last factors is trivial. The coefficient modules are 
then tensor products of the ${\cal S}_j$-modules for which the stability
result was proven in the last subsection , and so the stability result follows
for these also.  This completes the proof of Theorem \ahpi.10, and hence, the
proof of Theorem B. \hfill\za

\bigskip

\centerline{\bf Bibliography} 
\medskip 
\parskip = 0pt
\font\ninesl=cmsl9 
\font\bsc=cmcsc10 at 10truept 
\ninerm  

\item{[Akh]}{\bsc D.N. Akhiezer}, {\ninesl Homogeneous Complex Manifolds}, in
Enc. Math. Sci., Vol 10, Several Complex Variables IV, Springer-Verlag, 1990.
\item{[AJ]} {\bsc M.F. Atiyah, J.D. Jones}, {\ninesl Topological
aspects of Yang-Mills theory}, Comm. Math. Phys. 61 (1978), 
97-118. 
\item{[Bo1]}{\bsc A. Borel}, {\ninesl Linear Algebraic Groups, Second Enlarged
Edition}, Springer-Verlag, 1991.
\item{[Bo2]}{\bsc A. Borel}, {\ninesl Symmetric compact complex spaces}, Arch.
Math. 33 (1979) 49-56. (Reprinted in A. Borel Collected Papers, Vol III).
\item{[BHMM1]}{\bsc C.P. Boyer, J.C. Hurtubise, B.M. Mann, R.J. 
Milgram}, {\ninesl The topology of holomorphic maps into generalized flag
manifolds}, Acta Math. 173 (1994) 61-101.
\item{[BHMM2]}{\bsc C.P. Boyer, J.C. Hurtubise, B.M. Mann, R.J. 
Milgram}, {\ninesl The topology of instanton moduli spaces. I: The
Atiyah-Jones conjecture}, Ann. of Math. 137 (1993) 561-609.
\item{[Br]}{\bsc M. Brion}, {\ninesl Spherical varieties; an introduction}, in
Topological Methods in Algebraic Transformation Groups, H. Kraft, T. Petrie, G.
Schwarz, Editors, Birkh\"auser, Progress in Math. Vol 80, 1989.
\item{[BLV]}{\bsc M.  Brion, D. Luna, Th. Vust}, {\ninesl Espaces homog\`enes
sph\'eriques}, Invent.  Math. 84 (1986) 617-632.  
\item {[DP]} {\bsc C. DeConcini, C. Procesi}, {\ninesl Complete symmetric
varieties}, in Invariant Theory, Springer Lecture Notes 996, 1-45.
\item{[EW]}{\bsc J. Eells, J.C. Wood}, {\ninesl Maps of minimum energy}, J.
London Math. [2] 23 (1981) 303-310.
\item{[Ful1]}{\bsc W. Fulton}, {\ninesl Intersection Theory}, Springer-Verlag,
1984.
\item{[Ful2]}{\bsc W. Fulton}, {\ninesl Introduction to Toric Varieties},
Annals of Math. Studies, Princeton Univ. Press, 1993.
\item{[GH]}{\bsc P. Griffiths, J. Harris}, {\ninesl Principles of Algebraic
Geometry}, John Wiley \& Sons, New York, 1978.
\item{[Gra]} {\bsc J. Gravesen}, {\ninesl On the topology of spaces of
holomorphic maps}, Acta Math. 162 (1989), 247-286.
\item{[Gro]} {\bsc A. Grothendieck}, {\ninesl A General Theory of Fibre Spaces
with Structure Sheaf}, Lecture Notes, University of Kansas, 1958.
\item{[Gu1]} {\bsc M.A. Guest}, {\ninesl Topology of the space of 
absolute minima of the energy functional},  Amer. J. of Math.
106(1984), 21-42. 
\item{[Gu2]} {\bsc M.A.  Guest}, {\ninesl The topology of the space of rational
curves on a toric variety}, Acta Math. 174 (1995) 119-145.  
\item{[HM]} {\bsc J.C. Hurtubise, R.J. Milgram}, {\ninesl The Atiyah-Jones
conjecture for ruled surfaces}, J. fur die Reine und Ang. Math. 466, (1995),
111-143.
\item{[HO]} {\bsc A.T. Huckleberrry, E. Oeljeklaus}, {\ninesl Classification
theorems for almost homogeneous spaces}, Revue de l'Institut Elie Cartan,
Nancy, No. 9, 1984.
\item{[HR]} {\bsc P. Hilton, J. Roitberg} {\ninesl On the Zeeman comparison theorem for 
the homology of quasi-nilpotent fibrations}, Quarterly J. Math.  27 (1976), 433-444
\item{[Hu1]} {\bsc J.C.  Hurtubise}, {\ninesl Holomorphic maps of a riemann
surface into a flag manifold}, J. Diff. Geom. 43, (1996), 99-118.
\item{[Hu2]} {\bsc J.C.  Hurtubise}, {\ninesl Moduli spaces and particle
spaces}, in Gauge Theory and Symplectic Geometry,  Proc. of the
NATO ASI S\'eminaire de Math\'ematiques Sup\'erieures, 
J. Hurtubise and F.Lalonde eds., Kluwer, The Netherlands, 1997.
\item{[Hu3]} {\bsc J.C.  Hurtubise}, {\ninesl Stability theorems for moduli
spaces}, in CMS-SMC 50th Anniversary volume 3, Invited Papers, Ottawa, 1997,
153-172.

\item{[Huc]} {\bsc A.T. Huckleberry}, {\ninesl Actions of groups of holomorphic
transformations}, in Enc. Math. Sci., Vol 69, Several Complex Variables VI,
Springer-Verlag, 1990.  
\item{[Iit]} {\bsc S. Iitaka}, {\ninesl Algebraic Geometry}, Springer-Verlag,
New York, 1982.  
\item {[Jan]} {\bsc J.C Jantzen} {\ninesl Representations of Algebraic Groups}, Academic Press,
Orlando, 1987.
\item{[Ki1]} {\bsc F.C. Kirwan}, {\ninesl On spaces of maps from Riemann
surfaces to Grassmannians and applications to the cohomology of moduli of
vector bundles}, Ark. Math. 24(2) (1986), 221-275.  
\item{[Ki2]} {\bsc F.C. Kirwan}, {\ninesl Geometric invariant theory and the
Atiyah-Jones conjecture}, Sophus Lie Memorial Conference Proceedings, O.A.
Laudal and B. Jahren eds., Scandinavian University Press, Oslo, 1994.
\item{[Ko]} {\bsc J. Koll\'ar}, {\ninesl Rational Curves on Algebraic
Varieties}, Springer-Verlag, New York, 1996.
\item{[Leh]} {\bsc R. Lehmann}, {\ninesl Complex-symmetric spaces}, Ann. Inst.
Fourier, Grenoble 39(2) (1989), 373-416.
\item{[Mi]} {\bsc Y. Miyaoka}, {\ninesl Rational curves on algebraic
varieties}, Proc. Inter. Cong., 1994, Birkh\"auser, 680-689.
\item{[MM1]} {\bsc B.M. Mann, R.J. Milgram}, {\ninesl Some spaces of
holomorphic maps to complex Grassmann manifolds}, J. Diff. Geom. 33(2)
(1991), 301-324.  
\item{[MM2]} {\bsc B.M. Mann, R.J.  Milgram}, {\ninesl On the moduli space of
SU(n) monopoles and holomorphic maps to flag manifolds}, J. Diff. Geom. 38
(1993), 39-103.
\item{[Oel]} {\bsc E.  Oeljeklaus}, {\ninesl Fasthomogene
K\"ahlermannigfaltigleiten mit verschwindender ersten Bettizahl}, Manuscr.
Math. 7 (1972) 175-183. 
\item{[Se]} {\bsc G. Segal}, {\ninesl The topology of rational functions},
Acta Math. 143 (1979), 39-72.  
\item {[T]} {\bsc C.H. Taubes}, {\ninesl The stable topology of
self-dual moduli spaces}, J.  Diff.  Geom., 29 (1989), 163-230.
\item{[Ti1]} {\bsc Y.  Tian}, {\ninesl The based
$\scriptstyle{SU(n)}$-instanton moduli spaces}, Math.  Ann. 298 (1994),
117-139.  
\item{[Ti2]} {\bsc Y. Tian}, {\ninesl The Atiyah-Jones conjecture for the
classical Lie groups and Bott peroidicity}, J. Diff. Geom. 44 (1996), 178-199.

\ninerm 
\bigskip 
\bigskip 
 
\line{C.P. Boyer  \hfil  J. C. Hurtubise} 
\line{Department of Mathematics and Statistics \hfil   
Department of Mathematics} 
\line{University of New Mexico \hfil McGill University} 
\line{email: cboyer@math.unm.edu \hfil  
email: hurtubis@gauss.math.mcgill.ca} 
 
\medskip 
 
\line{R.J. Milgram \hfil} 
\line{Department of Mathematics \hfil}
\line{Stanford University \hfil} 
\line{email: milgram@gauss.stanford.edu \hfil} 

\bye